# Algebraic tensor products and internal homs of noncommutative $L_p$-spaces

Dmitri Pavlov
Mathematical Institute, University of Münster
pavlov@math.berkeley.edu

**Abstract.** We prove that the multiplication map $L_a(M) \otimes_M L_b(M) \to L_{a+b}(M)$ is an isometric isomorphism of (quasi)Banach $M$-$M$-bimodules. Here $L_a(M) = L^{1/a}(M)$ is the noncommutative $L_p$-space of an arbitrary von Neumann algebra $M$ and $\otimes_M$ denotes the algebraic tensor product over $M$ equipped with the (quasi)projective tensor norm, but without any kind of completion. Similarly, the left multiplication map $L_a(M) \to \operatorname{Hom}_M(L_b(M), L_{a+b}(M))$ is an isometric isomorphism of (quasi)Banach $M$-$M$-bimodules, where $\operatorname{Hom}_M$ denotes the algebraic internal hom without any continuity assumptions. In a forthcoming paper these results will be applied to $L_p$-modules introduced by Junge and Sherman, establishing explicit algebraic equivalences between the categories of right $L_p(M)$-modules for all $p \geq 0$.

## 1 Introduction.

1.1 In this paper $L_a$ means $L^{1/a}$, in particular $L_0 = L^\infty$ and $L_{1/2} = L^2$. See the next section for details.

1.2 The classical inequality discovered by Leonard Rogers in 1888 and improved by Otto Hölder in 1889 can be stated in modern language as follows. If $x \in L_a(M)$ and $y \in L_b(M)$ then $xy \in L_{a+b}(M)$ and $\|xy\| \leq \|x\| \cdot \|y\|$, where $M$ is an arbitrary von Neumann algebra and $a$ and $b$ belong to $\mathbf{C}_{\Re \geq 0}$. Rogers and Hölder proved this inequality for finite-dimensional commutative von Neumann algebras. The proof for arbitrary von Neumann algebras was written down only in 1984 by Hideki Kosaki [27, Lemma 1] using Uffe Haagerup's 1979 construction [15, Definition 1.7] of $L_a(M)$ for an arbitrary von Neumann algebra $M$. We suppress a lot of intermediate history here, in particular, the cases of commutative and semifinite von Neumann algebras, which were worked out in the intervening 96 years. An important subtlety here is that the space $L_a(M)$ only has a quasinorm for $\Re a > 1$. Quasinorms are mild generalizations of norms explained below, satisfying the quasitriangle inequality $\|x + y\| \leq c(\|x\| + \|y\|)$. For the space $L_a(M)$ ($\Re a \geq 1$) we can take $c = 2^{\Re a - 1}$.

1.3 Hölder's inequality can be expressed in terms of the two maps $L_a(M) \otimes L_b(M) \to L_{a+b}(M)$ ($x \otimes y \mapsto xy$) and $L_a(M) \to \operatorname{Hom}(L_b(M), L_{a+b}(M))$ ($x \mapsto (y \mapsto xy)$). The inequality then says that these maps are contractive, i.e., preserve or decrease the norm. The multiplication map is associative, hence the above maps factor through the maps $L_a(M) \otimes L_b(M) \to L_a(M) \otimes_M L_b(M)$ and $\operatorname{Hom}_M(L_b(M), L_{a+b}(M)) \to \operatorname{Hom}(L_b(M), L_{a+b}(M))$. We equip the space $L_a(M) \otimes_M L_b(M)$ with the factor (quasi)norm. The space $\operatorname{Hom}_M(L_b(M), L_{a+b}(M))$ is similarly equipped with the restriction (quasi)norm.

1.4 Since quasinorms and their friends appear whenever we consider $L_a$-spaces for $\Re a > 1$, we pause to explain some basic definitions of the theory of quasi-Banach spaces and their tensor products. See Kalton [42] for a nice survey of this area. If $a \geq 1$ is a real number, then an *a-seminorm* on a complex vector space $V$ is a function $x \in V \mapsto \|x\| \in \mathbf{R}_{\geq 0}$ such that $\|ax\| = |a| \cdot \|x\|$ for any $a \in \mathbf{C}$ and any $x \in V$, and $\|x + y\|^{1/a} \leq \|x\|^{1/a} + \|y\|^{1/a}$ for all $x$ and $y$ in $V$. An *a-norm* is a $a$-seminorm satisfying the usual nondegeneracy condition: for all $x \in V$ the relation $\|x\| = 0$ implies $x = 0$. In particular, $a = 1$ gives the usual definition of a norm. We use the letter $a$ instead of the traditional $p$ to remind the reader of the reparametrization $a = 1/p$, hence the term $a$-norm instead of the usual $p$-norm.

1.5 An *a-Banach space* is a complete $a$-normed (complex) vector space. If $a \in \mathbf{C}_{\Re \geq 0}$ and $\Re a \geq 1$, then $L_a(M)$ is an example of an $\Re a$-Banach space. For $\Re a \leq 1$ the space $L_a(M)$ is an ordinary Banach space.

1.6 A $p$-norm is an example of a quasinorm, which we define in the same way as a $p$-norm, but replace the inequality by $\|x + y\| \leq c(\|x\| + \|y\|)$, where $c \geq 1$ is some constant. In the case of a $p$-norm we can take $c = 2^{p-1}$. A theorem by Aoki and Rolewicz (see Section 2 in Kalton [42]) states that every quasinorm with constant $c$ is equivalent (with respect to the obvious notion of equivalence of quasinorms) to a $(1 + \log_2 c)$-norm.

1.7 The tensor product and the internal hom above turn out to be automatically complete. Even more is true: the maps $L_a(M) \otimes_M L_b(M) \to L_{a+b}(M)$ and $L_a(M) \to \operatorname{Hom}_M(L_b(M), L_{a+b}(M))$ are isometric isomorphisms of $M$-$M$-bimodules. This is the first main result of this paper.

1.8 **Theorem.** For any von Neumann algebra $M$ and for any $a$ and $b$ in $\mathbf{C}_{\Re\geq 0}$ the multiplication map $\mathrm{L}_a(M)\otimes_M \mathrm{L}_b(M) \to \mathrm{L}_{a+b}(M)$ and the left multiplication map $\mathrm{L}_a(M) \to \mathrm{Hom}_M(\mathrm{L}_b(M), \mathrm{L}_{a+b}(M))$ are isometric isomorphisms of (quasi)Banach $M$-$M$-bimodules. Here $\otimes_M$ denotes the algebraic tensor product (without any kind of completion) and $\mathrm{Hom}_M$ denotes the algebraic internal hom (without any kind of continuity assumption).

1.9 The underlying intuition behind these claims is that the spaces $\mathrm{L}_a(M)$ for any $a \in \mathbf{C}_{\Re\geq 0}$ and any von Neumann algebra $M$ can be thought of as (algebraically) cyclic right $M$-modules. As stated this is true only for $\Re a = 0$. For $\Re a > 0$ only a weaker statement is true: every finitely generated algebraic submodule of $\mathrm{L}_a(M)$ is cyclic. Countably generated submodules that are closed under countable sums are also cyclic. Also $\mathrm{L}_a(M)$ is *topologically* cyclic (i.e., it admits a dense cyclic submodule) whenever $M$ is $\sigma$-finite, but this is not important in what follows.

1.10 The original motivation for these results comes from a desire to construct a local two-dimensional functorial quantum field theory in the sense of Stolz and Teichner [44]. The relevant ideas, including the usage of $\mathrm{L}_a$-spaces, date back to Graeme Segal [45]. Such a theory, in particular, would send the point (an object in the two-dimensional bordism bicategory) to some von Neumann algebra $M$ (the hyperfinite type III factor is a popular choice) and the angle of measure $a \in \mathbf{R}_{\geq 0}$ (a 1-morphism in the bordism bicategory, which is the endomorphism of the point that rotates it by the angle $a$) to some $M$-$M$-bimodule $I_a$. Composition of such endomorphisms adds the corresponding angles, hence for such a field theory to exist we must have an isomorphism of $M$-$M$-bimodules $I_a \otimes_M I_b \to I_{a+b}$. As it turns out from the results of this paper, the spaces $\mathrm{L}_a(M)$ form such a family. One motivation for this particular choice of $I_a$ is that for $a = 2\pi$ we can glue together the two ends of the angle of measure $2\pi$ and obtain a 1-morphism $C$ (an endomorphism of the empty 0-bordism) that looks like an infinitesimal punctured disk, i.e., an infinitesimal circle. If we rescale $a$ by a factor of $2\pi$ and set $I_a = \mathrm{L}_{a/2\pi}(M)$, then $C$ is mapped to $\mathrm{L}_1(Z)$, where $Z$ is the center of $M$. There is a 2-morphism in the bordism bicategory from $C$ to the empty 1-bordism that looks like an infinitesimal (unpunctured) disk filling an infinitesimal punctured disk with the missing point. This 2-morphism yields a morphism of vector spaces from the vector space assigned to $C$ to the vector space of complex numbers, which in our case is given by the Haagerup trace $\mathrm{L}_1(Z) \to \mathbf{C}$.

## 2 Notation and terminology.

2.1 By a von Neumann algebra $M$ we mean a complex W*-algebra, i.e., a complex unital C*-algebra $M$ that admits a (necessarily unique) predual $M_*$. We do not assume that $M$ is represented on a Hilbert space. The symbol $\otimes_M$ denotes the algebraic tensor product of an algebraic right $M$-module and an algebraic left $M$-module over $M$ and $\otimes$ denotes $\otimes_{\mathbf{C}}$. Likewise, $\mathrm{Hom}_M$ denotes the algebraic internal hom of algebraic right $M$-modules and Hom denotes $\mathrm{Hom}_{\mathbf{C}}$. By a morphism of von Neumann algebras we mean a normal (i.e., ultraweakly continuous) unital *-homomorphism.

2.2 A weight $\mu$ on a von Neumann algebra $M$ is an ultraweakly lower semicontinuous (equivalently, normal) additive homogeneous map $\mu\colon M^+ \to [0, \infty]$; $\mu$ is *semifinite* if the set $\{x \in M^+ \mid \mu(x) \neq \infty\}$ generates $M$ and *faithful* if $\mu(p) = 0$ implies $p = 0$ for all $p \in M^+$. Here $M^+$ denotes the set of positive elements of $M$. In our definition weights and related objects like states, traces, and operator-valued weights are normal by definition. We denote by $\mathrm{W}_{\mathrm{sf}}(M)$ the set of all semifinite weights on a von Neumann algebra $M$.

2.3 Denote by $\mathbf{I}$ the set of imaginary complex numbers $\{p \in \mathbf{C} \mid \Re p = 0\}$ and by $\Re\colon \mathbf{C} \to \mathbf{R}$ and $\Im\colon \mathbf{C} \to \mathbf{I}$ the projections of $\mathbf{C} \cong \mathbf{R} \oplus \mathbf{I}$ onto $\mathbf{R}$ and $\mathbf{I}$ respectively. In our notation $\Im p \in \mathbf{I}$, i.e., the imaginary part is actually imaginary, not real. Denote by $\mathbf{R}^\times_{>0}$ the (multiplicative) Lie group of strictly positive real numbers and by $\mathbf{U}$ the (multiplicative) Lie group of complex numbers with absolute value 1. The projections of $\mathbf{C}^\times \cong \mathbf{R}^\times_{>0} \times \mathbf{U}$ onto $\mathbf{R}^\times_{>0}$ and $\mathbf{U}$ are denoted by $z \mapsto |z|$ and arg respectively. Finally, denote by $\mathbf{R}_{\geq 0}$ the (additive) Lie monoid of nonnegative real numbers and by $\mathbf{C}_{\Re\geq 0}$ the (additive) Lie monoid of complex numbers with a nonnegative real part.

2.4 In (noncommutative) integration there are three parameters that denote the same quantity. The first parameter is used in the theory of smooth manifolds to parametrize density bundles. We denote it by $a$. The second parameter is usually denoted by $p$ and was introduced by F. Riesz to parametrize $\mathrm{L}_p$-spaces. It is also used in the theory of quasi-Banach spaces. The third parameter is usually denoted by $t$ and is used in Tomita-Takesaki theory to parametrize modular automorphism groups and Radon-Nikodym cocycle derivatives. We have $p = 1/a$ and $t = -ia$. The choice of $p$ is highly unnatural because $L_p$-spaces do not form a **C**-graded algebra in this notation. Similarly, $t$ is unsuitable because the involution on $L_t$-spaces does not interact properly with the **C**-grading and the notation itself forces us to make a noncanonical choice of a square root of $-1$. The choice of $a$ does not suffer from these problems, hence we use it consistently throughout this paper. This convention forces us to make adjustments to the traditional notations mentioned above. We denote by $L_a$ what is usually denoted by $L^{1/a}$. Thus, $L_0$ and $L_{1/2}$ denote what is usually denoted by $L^\infty$ and $L^2$. We refer to the elements of $L_a$ as *a-densities*, to stress their connection with smooth densities in differential geometry. Likewise, modular automorphism groups and Connes' Radon-Nikodym cocycle derivatives (the relevant definitions will be reviewed below) are parametrized by the elements of **I**, not **R**. To remind the reader of this reparametrization we consistently avoid using the letter $p$ in this context, which results in a somewhat unconventional term $L_a$-space.

## 3 Smooth $L_a$-spaces.

3.1 In this section we review the smooth counterpart of the theory of commutative $L_a$-spaces. The cases $a = 0$ and $a = 1$ correspond to bounded smooth functions and integrable smooth 1-densities on a smooth manifold. Although this section is not formally necessary for the rest of the paper, it provides helpful intuition for the more sophisticated construction of $L_a$-spaces of noncommutative von Neumann algebras.

3.2 Suppose $V$ is a finite-dimensional complex vector space and $a \in \mathbf{C}$. By definition, $\mathrm{Dens}_a(V)$ is the one-dimensional complex vector space consisting of set-theoretical functions $x \colon \det(V)^\times \to \mathbf{C}$ such that $x(pg) = |p|^a x(g)$ for all $p \in \mathbf{C}^\times$ and $g \in \det(V)^\times$, where $\det(V)$ denotes the top exterior power of $V$ and $\det(V)^\times$ is the $\mathbf{C}^\times$-torsor of nonzero elements of $\det(V)$. The vector spaces $\mathrm{Dens}_a(V)$ can be organized into a **C**-graded unital *-algebra using pointwise multiplication and conjugation. We can already see the main isomorphisms in this case: the maps $\mathrm{Dens}_a(V) \otimes \mathrm{Dens}_b(V) \to \mathrm{Dens}_{a+b}(V)$ and $\mathrm{Dens}_a(V) \to \mathrm{Hom}(\mathrm{Dens}_b(V), \mathrm{Dens}_{a+b}(V))$ are injective linear maps of one-dimensional vector spaces, hence they are isomorphisms. Here $a$ and $b$ are arbitrary complex numbers.

3.3 The vector space $\mathrm{Dens}_0(V)$ is canonically isomorphic to **C** and the space $\mathrm{Dens}_1(V)$ is canonically isomorphic to the dual of $\det(V) \otimes \mathrm{or}(V)$, where $\mathrm{or}(V)$ denotes the complex orientation line of $V$, i.e., the vector space of all set-theoretical functions $x \colon \det(V)^\times \to \mathbf{C}$ such that $x(pg) = \arg(p)x(g)$ for all $p \in \mathbf{C}^\times$ and $g \in \det(V)^\times$.

3.4 Nonzero $a$-densities are invertible: if $x \in \mathrm{Dens}_a(V)$, then $x^{-1}$ is the unique element of $\mathrm{Dens}_{-a}(V)$ such that $xx^{-1} = 1 \in \mathrm{Dens}_0(V) \cong \mathbf{C}$. In fact, $x^{-1}$ can be constructed explicitly as the function $g \mapsto (x(g))^{-1}$.

3.5 If $\Re a = 0$ then the vector space $\mathrm{Dens}_a(V)$ has a canonical norm, which for an element $x \in \mathrm{Dens}_a(V)$ is given by the common absolute value of all elements in the image of $x$ considered as a function $x \colon \det(V)^\times \to \mathbf{C}$. The formula $x(pg) = |p|^a x(g)$ implies that all such elements have the same absolute value because $|p|^a \in \mathbf{U}$ for $a \in \mathbf{I}$.

3.6 If $\Im a = 0$ then it makes sense to talk about the positive part $\mathrm{Dens}_a^+(V)$ of $\mathrm{Dens}_a(V)$, which consists of all functions in $\mathrm{Dens}_a(V)$ with values in $\mathbf{R}_{\geq 0}$. Thus $\mathrm{Dens}_a(V)$ is canonically oriented for all $a \in \mathbf{R}$. For every $b \in \mathbf{C}_{\Re \geq 0}$ we have the power map $\mathrm{Dens}_a^+(V) \to \mathrm{Dens}_{ab}(V)$, which lands in the positive part of $\mathrm{Dens}_{ab}(V)$ whenever $\Im b = 0$. Explicitly, if $x \in \mathrm{Dens}_a^+(V)$, then $x^b \in \mathrm{Dens}_{ab}(V)$ is defined by the formula $x^b(g) := (x(g))^b$. Here $0^b := 0$ even for $\Re b = 0$, thus for $x = 0$ we have $x^b = 0$ for all $b \in \mathbf{C}_{\Re \geq 0}$. The above power map can also be defined for $\Re b < 0$, but we must restrict its domain to the nonzero elements of $\mathrm{Dens}_a^+(V)$.

3.7 The power map can be used to define the absolute value and polar decomposition for densities. Given $x \in \mathrm{Dens}_a(V)$ for some $a \in \mathbf{C}$, we define $|x| := (xx^*)^{1/2} \in \mathrm{Dens}_{\Re a}^+(V)$ and $\arg(x) := x/|x| \in \mathrm{Dens}_{\Im a}(V)$. The argument function is only defined if $x$ is invertible, i.e., nonzero. In this case $\arg(x)$ is unitary: $\arg(x)\arg(x)^* = 1 \in \mathrm{Dens}_0(V) \cong \mathbf{C}$.

3.8 The power map also allows us to represent any $a$-density $x$ as $x = pg^a$ using a 0-density $p$ (i.e., a complex number) and a positive 1-density $g$. Indeed, if $\Re a > 0$ and $x \in \mathrm{Dens}_a(V)$ is a nonzero density, set $g := |x|^{1/\Re a} \in \mathrm{Dens}_1^+(V)$ and $p := \arg(x)g^{-\Im a} \in \mathrm{Dens}_0(V) \cong \mathbf{C}$ and observe that $pg^a = \arg(x)g^{-\Im a}g^a = \arg(x)g^{\Re a} = \arg(x)|x| = x$. For $\Re a = 0$ and $x \in \mathrm{Dens}_a(V)$ a nonzero density we have $x = |x| \cdot \arg(x)$, where $|x| \in \mathrm{Dens}_0(V) \cong \mathbf{C}$ and the unitary element $\arg(x) \in \mathrm{Dens}_a(V)$ can in its turn be represented as $\arg(x) = g^a$ for some nonunique noncanonical $g \in \mathrm{Dens}_1^+(V)$. In fact, the set of all such $g$ is a torsor over the kernel $\mathbf{P}$ of the map $\exp\colon \mathbf{I} \to \mathbf{U}$, which can be canonically identified with "imaginary integers" (elements of $\mathbf{I}$ with an integer absolute value) by dividing by $2\pi$. Indeed, if $a \in \mathbf{I}^\times$ and $u \in \mathrm{Dens}_a(V)$ is unitary, then for any $g \in \mathrm{Dens}_1^+(V)$ such that $g^a = u$ the element $p \cdot g := \exp(pa^{-1})g \in \mathrm{Dens}_1^+(V)$ again satisfies $(p \cdot g)^a = u$ for any $p \in \mathbf{P}$.

3.9 The above constructions can be done in smooth families. In particular, if we apply them to the complex tangent bundle of a smooth manifold $X$ we obtain the line bundles of densities $\mathrm{Dens}_a(X)$ for every $a \in \mathbf{C}$.

3.10 The power map, the argument, and the absolute value of a smooth density are only defined on the open subset of $X$ on which the density is nonvanishing because the map $t \in [0, \infty) \mapsto t^a \in \mathbf{C}$ is always smooth on $(0, \infty)$, but it is only smooth at 0 if $a$ is a nonnegative integer. We could extend the result of these maps by zero to the entire manifold $X$, similarly to how polar decomposition of operators is defined, but then we would be forced to consider discontinuous sections, which can be done by passing to the space of measurable $a$-densities, as explained below.

3.11 Similarly, the above decomposition of a nonvanishing $a$-density $x$ in the form $x = pg^a$ works as expected for $\Re a > 0$ because all constructions are canonical, but for $\Re a = 0$ we have a nontrivial topological obstruction for representing a unitary $a$-density as a positive 1-density raised to the power $a$. Indeed, given a unitary element $x \in \mathrm{C}^\infty(\mathrm{Dens}_a(X))$ we can construct a principal $\mathbf{P}$-bundle whose fiber at any given point is the $\mathbf{P}$-torsor of all solutions $g$ to the equation $g^a = x$ at the given point, as described above. A global solution to the equation $g^a = x$ exists if and only if this bundle is trivial.

3.12 Alternatively, one can say that unitary $a$-densities are smooth sections of the principal $\mathbf{U}$-bundle whose fibers are unitary $a$-densities at a given point and the $\mathbf{U}$-action is given by multiplication. Likewise, positive 1-densities can be thought of (if $a \in \mathbf{I}^\times$) as smooth sections of the principal $\mathbf{I}$-bundle whose fibers are positive 1-densities at a given point and the action by an element $p \in \mathbf{I}$ is given by the multiplication by $\exp(pa^{-1})$. The power map induces a morphism from the latter bundle to the former associated to the morphism $\exp\colon \mathbf{I} \to \mathbf{U}$ of the corresponding structure groups, hence we have an induced map on global sections and on their (smooth) homotopy classes. However, there is just one homotopy class of sections of a principal $\mathbf{I}$-bundle, whereas a principal $\mathbf{U}$-bundle can have many such classes, and exactly one of these classes (namely, the image of the only homotopy class of 1-densities) allows for global solutions to the equation $g^a = x$. Again the lifting problem disappears once we pass to measurable $a$-densities.

3.13 For any smooth manifold $X$ we have the integration map $\int\colon \mathrm{C}_{\mathrm{cs}}^\infty(\mathrm{Dens}_1(X)) \to \mathbf{C}$ ($\mathrm{C}_{\mathrm{cs}}^\infty$ denotes compactly supported smooth sections), which can be defined as the composition of the canonical factor map $\mathrm{C}_{\mathrm{cs}}^\infty(\mathrm{Dens}_1(X)) \to \mathrm{H}_{\mathrm{cs}}^{\mathrm{top}}(X, \mathrm{or}(X))$ for the twisted de Rham cohomology with coefficients in $\mathrm{or}(X)$, the Poincaré duality isomorphism $\mathrm{H}_{\mathrm{cs}}^{\mathrm{top}}(X, \mathrm{or}(X)) \to \mathrm{H}_0(X, \mathbf{C})$, and the pushforward map in ordinary homology $\mathrm{H}_0(X, \mathbf{C}) \to \mathrm{H}_0(\mathrm{pt}, \mathbf{C}) = \mathbf{C}$ induced by the unique map $X \to \mathrm{pt}$. The integration map preserves positivity: $\int \mathrm{C}_{\mathrm{cs}}^\infty(\mathrm{Dens}_1^+(X)) \subset \mathbf{R}_{\geq 0}$.

3.14 We extend the integration map by continuity to the space of *integrable* smooth 1-densities, which are defined as those smooth 1-densities $x$ for which $\int$ is bounded on the set of all compactly supported 1-densities $y$ such that $0 \leq y \leq |x|$. An $a$-density $x$ for $\Re a > 0$ is *integrable* if $|x|^{1/\Re a}$ is an integrable 1-density (here and below we restrict to the open subset of $X$ on which $x$ is nonvanishing). The *measurable topology* on integrable $a$-densities is induced by the norm (for $\Re a \leq 1$) or quasinorm (for $\Re a > 1$) defined by the formula $\|x\| = \left(\int |x|^{1/\Re a}\right)^{\Re a} \in \mathbf{R}_{\geq 0}$. An $a$-density $x$ for $\Re a = 0$ is *bounded* if its fiberwise norm is a bounded function. The *measurable topology* on bounded $a$-densities is the weak topology induced by integrable smooth $(1-a)$-densities via the pairing defined by the formula $(x, y) = \int(\bar{x}y)$. We adopt a common notation for both of these cases and denote the space of bounded (for $\Re a = 0$) or integrable (for $\Re a > 0$) smooth $a$-densities equipped with the measurable topology by $\mathrm{C}_{\mathrm{L}}^\infty(\mathrm{Dens}_a(X))$. The subscript L refers to the fact that these spaces should be thought of as smooth versions of $\mathrm{L}_a$-spaces. There is no reasonable topology on $a$-densities for $\Re a < 0$, which explains why the spaces $\mathrm{L}_a$ are considered only for $\Re a \geq 0$. In particular, the above formula for the (quasi)norm does not work for $\Re a < 0$ because the power map is defined only for nonvanishing densities if $\Re a < 0$.

3.15 For any $a \in \mathbf{C}_{\Re \geq 0}$ we define the topological vector space $\mathrm{L}_a(X)$ as the completion of $\mathrm{C}_{\mathrm{L}}^\infty(\mathrm{Dens}_a(X))$. Both the multiplication map $\mathrm{C}_{\mathrm{L}}^\infty(\mathrm{Dens}_a(X)) \otimes_{\mathrm{C}_{\mathrm{L}}^\infty(X)} \mathrm{C}_{\mathrm{L}}^\infty(\mathrm{Dens}_b(X)) \to \mathrm{C}_{\mathrm{L}}^\infty(\mathrm{Dens}_{a+b}(X))$ and the left multiplication map $\mathrm{C}_{\mathrm{L}}^\infty(\mathrm{Dens}_a(X)) \to \mathrm{Hom}_{\mathrm{C}_{\mathrm{L}}^\infty(X)}(\mathrm{C}_{\mathrm{L}}^\infty(\mathrm{Dens}_b(X)), \mathrm{C}_{\mathrm{L}}^\infty(\mathrm{Dens}_{a+b}(X)))$ can be seen to be continuous isomorphisms and therefore can be extended to the corresponding completions in the measurable topology. This fact proves the main theorems for all separable commutative von Neumann algebras and their products because every coproduct of separable measurable spaces is the underlying measurable space of some smooth manifold. Smooth manifolds are assumed to be paracompact and Hausdorff but not necessarily second countable. (A paracompact Hausdorff smooth manifold is second countable if and only if it has at most countably many connected components; we need smooth manifolds with arbitrary many connected components.)

### 4 History of noncommutative integration.

4.1 In this section we briefly review the development of noncommutative integration and indicate how our approach to $L_a$-spaces fits into the general picture.

4.2 In 1910 Friedrich Riesz defined $L_a$ for the case of ordinary (commutative) measurable spaces. In 1929 John von Neumann introduced noncommutative measurable spaces (von Neumann algebras), which raised the problem of extending Riesz' definition to the noncommutative case. In 1953 Irving Segal [2, Definitions 3.3 and 3.7] defined noncommutative $L_1$ and $L_{1/2}$ spaces. Ogasawara and Yoshinaga [4, 5] explored some of their properties and Stinespring [7] developed the theory further. Kunze [6, Definition 3.1] defined $L_a$ for all $a \in [0,1]$, extending Segal's definition. Yeadon [13] also explored their properties in a similar framework, proving Hölder's inequality and duality, apparently independently from Kunze, given the absence of references to his paper. Dixmier [3, Section 3] published his own definition of $L_a$ for $a \in [0,1]$ with an accompanying theory simultaneously with Segal, defining the spaces $L_a$ as abstract completions. Nelson [12] wrote a nice exposition of the theory using the measure topology introduced by Stinespring. His approach is similar to Dixmier's even though he doesn't cite any of Dixmier's papers. These authors assumed that the von Neumann algebra under consideration is semifinite. The existence of a faithful semifinite trace $\tau$ is crucial for these constructions because elements in $L_a$ are essentially represented as multiples of $\tau^a$ and the fact that $\tau$ is a trace as opposed to a weight is used heavily. Haagerup [15, Definition 1.7] was the first one to define $L_a$ for $a \in [0,1]$ for an arbitrary von Neumann algebra, using crossed products and the modular flow. Terp [22] wrote up an elaborate exposition of his construction and the related theory. Meanwhile Connes [19, the paragraph after the proof of Corollary 18] gave his own definition using spatial derivatives, which are closely related to the relative modular flow on invertible bimodules, as explained by Yamagami in [35, Subsection "Spatial Derivatives" in Section "Modular Extension"]. Hilsum [21] wrote up the corresponding theory. Araki and Masuda [25] gave a definition for the $\sigma$-finite case based on the relative modular theory and the standard representation of a von Neumann algebra, which was soon extended by Masuda [26] to the general case.

4.3 Another series of approaches uses Calderón's complex interpolation method. Kosaki [28] gave the first definition, which was restricted to the case of $\sigma$-finite von Neumann algebras. Terp [23] extended Kosaki's construction to all von Neumann algebras. Izumi [36, 37] defined for every $a \in (0,1)$ a family of spaces parametrized by $\alpha \in \mathbf{C}$ that are functorially isomorphic to each other and to $L_a(M)$ in such a way that the constructions of Kosaki and Terp correspond to $\alpha \in [-1/2, 1/2]$ and $\alpha = 0$. Leinert [31, 32] developed another interpolation-based approach by defining a noncommutative analog of the upper integral in the semifinite case and later [33] extended it to the general case. Pisier and Xu [40] wrote a survey summarizing the results mentioned above. Unfortunately, all these approaches only consider $a$-densities for $a \in [0,1]$ and not for $a \in \mathbf{C}_{\Re \geq 0}$, which precludes their use in our paper.

4.4 All of the above definitions depend on a choice of a faithful semifinite weight $\mu$ on $M$. However, the spaces $L_a(M, \mu)$ for different $\mu$ are canonically isomorphic to each other. Thus we can define $L_a(M)$ as the limit (or the colimit, because all maps are isomorphisms) of $L_a(M, \mu)$ for all faithful semifinite weights $\mu$ on $M$. See the next section for the details of this construction in the commutative case. Kosaki [20] gave a weight-independent definition of $L_a(M)$ for $a \in (0,1]$ by introducing a new additive structure on the predual of $M$ using the relative modular theory. Earlier Haagerup [14] gave a weight-independent characterization of $L_{1/2}(M)$ (also known as the standard form of $M$) for an arbitrary von Neumann algebra $M$.

4.5 Yamagami [34] reformulated the original results of Haagerup [15] in a more convenient algebraic setting of modular algebras and defined $L_a(M)$ for an arbitrary $a \in \mathbf{C}_{\Re \geq 0}$. Sherman [39] along with Falcone and Takesaki [38] give alternative expositions. The principal idea of the modular algebra approach is to construct a **C**-graded *-algebra (the modular algebra of $M$) whose component in grading $a \in \mathbf{C}$ happens to be the space $\mathrm{L}_a(M)$ (zero for $\Re a < 0$). We use the language of modular algebras systematically throughout this paper. We note in passing that allowing $a$ to be purely imaginary (as opposed to forcing $a$ to be real) is crucial for establishing relations with the Tomita-Takesaki theory, as explained below.

4.6 In terms of modular algebras (see below for the relevant details) Kosaki's definition amounts to reinterpreting an element $x = py \in M_* \cong \mathrm{L}_1(M)$ as $py^a \in \mathrm{L}_a(M)$ and introducing new algebraic structures on $M_*$ accordingly, obtaining $\mathrm{L}_a(M)$. Here $p \in M$ and $y \in M_*^+$ form the right polar decomposition of $x = py$. Falcone and Takesaki's approach amounts to constructing $\mathrm{L}_a(M)$ for $a \in \mathbf{I}$ as the set of equivalence classes of pairs $(p, y)$, where $p \in M$, $y \in \mathrm{W}_{\mathrm{sf}}(M)$, and $(p, y) \sim (q, z)$ whenever $py^a = qz^a$.

## 5 Commutative $\mathrm{L}_a$-spaces.

5.1 In this section we establish a connection with the classical theory of $\mathrm{L}_a$-spaces. Unfortunately, the underlying notions of measure theory are not spelled out explicitly in the literature, even though they are crucial for the Gelfand duality for commutative von Neumann algebras, which itself seems to be lacking a decent exposition. In particular, the notion of a *measurable* space (a set equipped with a $\sigma$-algebra of measurable sets and a $\sigma$-ideal of negligible sets; also known as a *measure class*), as opposed to a *measure* space (a set equipped with a $\sigma$-algebra of measurable sets) and to a *measured* space (a set equipped with a $\sigma$-algebra of measurable sets and a measure), as well as the notion of a morphism of measurable spaces (an equivalence class of measurable functions modulo equality on the complement of a negligible set) are rarely mentioned in expository texts on measure theory. Thus we have to review the relevant definitions, making this section much longer than it could have been.

5.2 A *measurable space* is a triple $(X, M, N)$, where $X$ is a set, $M$ is a $\sigma$-algebra of *measurable subsets* of $X$, and $N \subset M$ is a $\sigma$-ideal of *negligible sets*. For the sake of simplicity we assume that our measurable spaces are *complete*, i.e., every subset of an element of $N$ is again an element of $N$. The inclusion functor from the category of complete measurable spaces to the category of measurable spaces is an equivalence, thus we do not lose anything by restricting ourselves to complete measurable spaces.

5.3 If $(X, M, N)$ and $(Y, P, Q)$ are measurable spaces, then a map $f \colon X \to Y$ is *measurable* if the preimage of every element of $P$ is an element of $M$. Measurable maps are closed under composition. A measurable map is *nonsingular* if the preimage of every element of $Q$ is an element of $N$. Nonsingular measurable maps are also closed under composition. Two measurable maps $f$ and $g$ are *equivalent* if $\{x \in X \mid f(x) \neq g(x)\} \in N$. Composition of nonsingular measurable maps preserves this equivalence relation. We define a *morphism of measurable spaces* as an equivalence class of nonsingular measurable maps.

5.4 A *measure* on a measurable space $(X, M, N)$ is a $\sigma$-additive map $\mu \colon M \to [0, \infty]$ such that $\mu(A) = 0$ for all $A \in N$. A measure $\mu$ is $\sigma$-*finite* if the union of all $A \in M$ such that $\mu(A) \neq \infty$ equals $X$. A measure $\mu$ is *faithful* if $\mu(A) = 0$ implies $A \in N$ for all $A \in M$. Abusing the language we say that a (complex-valued) *finite measure* on a measurable space $(X, M, N)$ is a $\sigma$-additive map $\mu \colon M \to \mathbf{C}$ such that $\mu(A) = 0$ for all $A \in N$. Denote the complex vector space of all finite measures on a measurable space $Z$ by $\mathrm{L}_1(Z)$. This is a Banach space via the norm defined below.

5.5 A (complex-valued) *function* on a measurable space $Z$ is an equivalence class of (possibly singular) measurable maps from $Z$ to the measurable space of complex numbers equipped with the $\sigma$-algebra of Lebesgue measurable sets and the standard $\sigma$-ideal of negligible sets. We remark that a function on a measurable space $Z$ can in fact be defined as a morphism from $Z$ to a certain measurable space, which contains the measurable space of complex numbers mentioned above and a lot of other stuff, e.g., an isolated point for every complex number, which allows us to get constant functions as morphisms of measurable spaces. Our use of the above ad hoc definition with singular maps is motivated exclusively by our desire to avoid developing here the theory that constructs this measurable space.

5.6 A function is *bounded* if at least one representative of its equivalence class is bounded. Denote the set of all bounded functions on a measurable space $Z$ by $\mathrm{L}_0(Z)$, which is a C*-algebra in a natural way. Every element $\mu$ of $\mathrm{L}_1(Z)$ yields a unique norm-continuous linear functional on $\mathrm{L}_0(Z)$, whose value on the characteristic function of a set $A \in M$ is $\mu(A)$. Now we equip $\mathrm{L}_1(Z)$ with the norm $\|\mu\| := \sup_{f \in \mathrm{L}_0(Z)} \mu(f)$, where $\|f\| \leq 1$.

5.7 Irving Segal in [1, Theorem 5.1] and Kelley in [8] proved that for a measurable space $Z = (X, M, N)$ the following properties are equivalent.

- The Boolean algebra $M/N$ of equivalence classes of measurable sets is complete, i.e., admits arbitrary suprema.
- The lattice of all real functions on $Z$ is Dedekind complete. An ordered set is *Dedekind complete* if any of its nonempty subsets bounded from above (respectively below) has a supremum (respectively infimum).
- The lattice of all bounded real functions on $Z$ is Dedekind complete.
- $Z$ has the Radon-Nikodym property: for any two faithful measures $\mu$ and $\nu$ on $Z$ we have $\mu = (\mathrm{D}\mu : \mathrm{D}\nu)\nu$ for a (unique unbounded strictly positive) function $(\mathrm{D}\mu : \mathrm{D}\nu)$ on $Z$.
- $Z$ has the Riesz representation property: the functorial evaluation map from $\mathrm{L}_0(Z)$ to the dual of the Banach space $\mathrm{L}_1(Z)$ is an isomorphism.
- $Z$ has the Hahn decomposition property: for every element $\mu \in \mathrm{L}_1(Z)$ we can represent $Z$ in the form $Z = Z_0 \sqcup Z_1$, where $\mu$ is positive on $Z_0$ and $-\mu$ is positive on $Z_1$.
- $\mathrm{L}_0(Z)$ is a von Neumann algebra.

5.8 If $Z$ satisfies any of these properties then we say that $Z$ is *localizable*. In particular, every $\sigma$-finite measurable space (i.e., a measurable space that admits a faithful finite measure) is localizable, but not vice versa. However, any localizable measurable space can be (noncanonically) represented as a coproduct of $\sigma$-finite spaces. Thus Segal's theorem shows that one has to work with localizable measurable spaces to get the usual theorems of measure theory. Moreover, the category of localizable measurable spaces is contravariantly equivalent to the category of commutative von Neumann algebras. The equivalence functor sends a measurable space $Z$ to $\mathrm{L}_0(Z)$ and a morphism of measurable spaces to the corresponding pullback map for functions. This theorem is essentially the von Neumann algebra analog of the Gelfand duality for commutative C*-algebras. Moreover, the dual of $\mathrm{L}_1(Z)$ is canonically isomorphic to $\mathrm{L}_0(Z)$ and the dual of $\mathrm{L}_0(Z)$ in the weak topology induced by $\mathrm{L}_1(Z)$ (the ultraweak topology) is canonically isomorphic to $\mathrm{L}_1(Z)$. Measures on $Z$ are canonically identified with weights on $\mathrm{L}_0(Z)$ (normal $[0, \infty]$-valued functionals on $\mathrm{L}_0^+(Z)$). Henceforth we include the property of localizability in the definition of a measurable space.

5.9 If $\mu$ is a faithful measure on a measurable space $Z$ (every measurable space admits such a measure, which need not be finite), define $\mathrm{L}_1(Z, \mu)$ as the space of all functions $f$ on $Z$ such that $\mu(|f|)$ is finite, where $|f| := (ff^*)^{1/2}$ is the absolute value of $f$. We introduce a Banach norm on $\mathrm{L}_1(Z, \mu)$ via $\|f\| := \mu(|f|)$. Then $\mathrm{L}_1(Z, \mu)$ and the space of finite measures $\mathrm{L}_1(Z)$ are functorially isometrically isomorphic via the multiplication map $f \in \mathrm{L}_1(Z, \mu) \mapsto f\mu \in \mathrm{L}_1(Z)$, in particular $\mathrm{L}_1(Z, \mu)$ is canonically isometrically isomorphic to $\mathrm{L}_1(Z, \nu)$ via the multiplication by the Radon-Nikodym derivative $(\mathrm{D}\mu : \mathrm{D}\nu)$.

5.10 For an arbitrary $a \in \mathbf{C}_{\Re \geq 0}$ we define $\mathrm{L}_a(Z, \mu)$ as the set of all functions $f$ on $Z$ such that $\left(\mu(|f|^{1/\Re a})\right)^{\Re a}$ is finite if $\Re a > 0$ or $\sup |f|$ is finite if $\Re a = 0$. The above expressions define a norm on $\mathrm{L}_a(Z, \mu)$ for $\Re a \leq 1$ and a quasinorm for $\Re a > 1$. If $\mu$ and $\nu$ are two faithful measures on $Z$, then we have a canonical isometric isomorphism $f \in \mathrm{L}_a(Z, \mu) \mapsto f(\mathrm{D}\mu : \mathrm{D}\nu)^a \in \mathrm{L}_a(Z, \nu)$. The space of $a$-densities $\mathrm{L}_a(Z)$ is the limit (or the colimit, because all maps are isomorphisms) of $\mathrm{L}_a(Z, \mu)$ for all $\mu$. The individual spaces $\mathrm{L}_a(Z, \mu)$ do not depend on the imaginary part of $a$, but the isomorphisms between them do, hence $\mathrm{L}_a(Z)$ is noncanonically isomorphic to $\mathrm{L}_{\Re a}(Z)$, and choosing such an isomorphism is equivalent to choosing a measure on $Z$. The spaces $\mathrm{L}_a(Z)$ and $\mathrm{L}_{\Re a}(Z)$ are no longer isomorphic as $\mathrm{L}_0(Z)$-$\mathrm{L}_0(Z)$-bimodules in the noncommutative case.

# 6 Noncommutative $L_a$-spaces.

6.1 In this section we define the Takesaki dual functor $M \mapsto \hat{M}$ from the category of von Neumann algebras and faithful semifinite operator-valued weights (which should be thought of as noncommutative analogs of maps equipped with fiberwise measures used to define pushforwards, see below for a precise definition) to the category of **I**-graded von Neumann algebras (see below for the definition of gradings). We then define $L_a(M)$ for $a \in \mathbf{I}$ as the $a$-graded component of $\hat{M}$, which gives a functor on the same category. Furthermore, $\hat{M}$ has a canonical faithful semifinite trace $\tau$, which we use to complete $\hat{M}$ in the $\tau$-measurable topology and obtain a **C**-graded extended von Neumann algebra (a slightly relaxed version of a von Neumann algebra that allows unbounded elements) $\bar{M}$, whose $a$-graded component is the space $L_a(M)$ for $a \in \mathbf{C}_{\Re \geq 0}$ and zero for all other $a$. These spaces keep the usual properties of their commutative versions, in particular, they have (quasi)norms that turn them into (quasi)Banach spaces and their elements have left and right polar decompositions and supports.

6.2 For the sake of intuition, we pause briefly to explain (in a somewhat informal way) the smooth commutative analog of this construction. We send the commutative *-algebra of bounded smooth functions on a smooth manifold $M$ to the commutative *-algebra $\hat{M}$ of distributional sections of the line bundle of imaginary densities over $M \times \mathbf{I}$ whose pullback to $M \times a$ ($a \in \mathbf{I}$) is the line bundle of $a$-densities over $M$. Informally, sections have to be smooth and bounded in the direction of $M$ and should have smooth bounded Fourier transform in the direction of $\mathbf{I}$. The unit is given by the delta distribution along $M \times 0$. The multiplication is given by the convolution in the direction of $\mathbf{I}$ and the multiplication in the direction of $M$. The involution is defined by precomposing with the inverse map in the direction of $\mathbf{I}$ and conjugating in the direction of $M$. The component in grading $a \in \mathbf{I}$ consists of sections supported on $M \times a$, which are smooth bounded sections of the line bundle of $a$-densities of $M$, i.e., the smooth $L_a$-space of $M$.

6.3 Since continuous gradings on von Neumann algebras are essential in the discussion below, we discuss them first.

6.4 **Definition.** If $G$ is an abelian locally compact topological group, then a $G$-*grading* on a von Neumann algebra $M$ is a morphism of groups $\theta: \hat{G} \to \mathrm{Aut}(M)$ such that the map $\hat{g} \in \hat{G} \mapsto \theta_{\hat{g}}(p) \in M$ is ultraweakly continuous for all $p \in M$. Here $\hat{G} := \mathrm{Hom}(G, \mathbf{U})$ is the Pontrjagin dual group of $G$ equipped with the compact-open topology (recall that $\mathbf{U}$ is the group of unitary complex numbers). A *morphism* of $G$-graded von Neumann algebras is a morphism of von Neumann algebras that commutes with $\theta$. For $g \in G$ we define the $g$-*graded component* of $M$ as the set of all elements $p \in M$ such that $\theta_{\hat{g}}(p) = \hat{g}(g)p$ for all $\hat{g} \in \hat{G}$. This construction extends to a functor: a morphism of $G$-graded von Neumann algebras $M \to N$ induces a morphism of their $g$-graded components because it commutes with the gradings.

6.5 In our case $G = \mathbf{I}$ and $\hat{G} = \mathrm{Hom}(\mathbf{I}, \mathbf{U})$ is identified with $\mathbf{R}$ via the following map: $\hat{g} \in \mathbf{R} \mapsto (g \in \mathbf{I} \mapsto \exp(-\hat{g}g) \in \mathbf{U}) \in \mathrm{Hom}(\mathbf{I}, \mathbf{U})$. The **I**-grading on the algebra $\hat{M}$ constructed below is also known as the *scaling automorphism group* or the *noncommutative flow of weights*. In terms of the above smooth manifold construction, acting by $\theta_{\hat{g}}$ on some section $p$ of the line bundle of imaginary densities over $M \times \mathbf{I}$ amounts to multiplying $p$ by the composition of the projection $M \times \mathbf{I} \to \mathbf{I}$ and $\hat{g}: \mathbf{I} \to \mathbf{U}$,

6.6 We can informally describe the core $\hat{M}$ of a von Neumann algebra $M$ as the **I**-graded von Neumann algebra $\hat{M}$ generated by $M$ in grading 0 (in particular, we talk about elements of $M$ as if they are also elements of $\hat{M}$) and symbols $\mu^a$ in grading $a$, where $a \in \mathbf{I}$ and $\mu \in \mathrm{W}_{\mathrm{sf}}(M)$, subject to the following relations.

- For all $\mu \in \mathrm{W}_{\mathrm{sf}}(M)$ the map $a \in \mathbf{I} \mapsto \mu^a \in \mathbf{U}(p\hat{M}p)$ is a continuous morphism of groups. As we explained before, $\mathrm{W}_{\mathrm{sf}}(M)$ denotes the set of semifinite weights on $M$. In the above formula $p$ denotes the support of $\mu$ (i.e., the minimum projection $p$ such that $\mu(pxp) = \mu(x)$ for all $x \in M^+$) and $\mathbf{U}(N)$ denotes the topological group of unitary elements of a von Neumann algebra $N$ equipped with the ultraweak topology.
- For all $\mu \in \mathrm{W}_{\mathrm{sf}}(M)$, $p \in M$, and $a \in \mathbf{I}$ we have $\mu^a p \mu^{-a} = \sigma_a^\mu(p)$, where $\sigma^\mu$ denotes the modular automorphism group of $\mu$ (see Definition VIII.1.3 in Takesaki [41]).
- For all $\mu \in \mathrm{W}_{\mathrm{sf}}(M)$, $\nu \in \mathrm{W}_{\mathrm{sf}}(M)$, and $a \in \mathbf{I}$ we have $\mu^a \nu^{-a} = (\mathrm{D}\mu : \mathrm{D}\nu)_a$, where $(\mathrm{D}\mu : \mathrm{D}\nu)$ denotes Connes' Radon-Nikodym cocycle derivative of $\mu$ with respect to $\nu$ (see Definition VIII.3.20 in Takesaki [41], which is easily extended to the case of nonfaithful $\nu$).

6.7 This construction can be readily interpreted in terms of smooth manifolds. We observe that all weights in the construction can be forced to be faithful without altering the resulting algebra $\hat{M}$. Smooth bounded functions on $M$ embed into the corresponding algebra $\hat{M}$ of sections of the line bundle of imaginary densities over $M \times \mathbf{I}$ as sections supported at $M \times 0$. Likewise, $\mu^a$ is a unitary section that depends continuously on $a$ and is supported at $M \times a$ for any $a \in \mathbf{I}$ and any nonvanishing 1-density $\mu$. Similarly, $\mu^a \nu^{-a}$ is a section supported at $M \times 0$ that corresponds to the $a$th power of the smooth function on $M$ defined as the ratio $\mu/\nu$, i.e., $\mu = (\mu/\nu)\nu$. However, due to the commutativity of the algebras involved we have $\mu^a p \mu^{-a} = \mu^a \mu^{-a} p = p$, hence $\sigma_a^\mu$ is always the identity automorphism. We obtain nontrivial smooth examples of $\sigma^\mu$ once we move into the world of higher smooth stacks, for example, Lie groupoids or even just Lie groups. In these cases the automorphisms $\sigma_a^\mu$ can be computed using *modular functions* (from which the names like modular algebra and modular automorphism group are derived). Unfortunately, considering these examples in any detail would lead us too far astray.

6.8 The last two relations in the definition of the core of $M$ could in fact be turned into *definitions* of $\sigma^\mu$ and $(\mathrm{D}\mu : \mathrm{D}\nu)$ once we have constructed the modular algebra. Of course, most traditional approaches to modular algebras require Tomita-Takesaki modular automorphisms and Connes' Radon-Nikodym derivatives to be constructed first, resulting in a circular dependence. However, as explained in a forthcoming paper by the author [47], it is possible to construct the **C**-graded extended von Neumann algebra $\bar{M}$ directly in a single step as the free **C**-graded extended von Neumann algebra generated by $M$ in degree 0 and $M_*$ in degree 1. In particular, this construction does not use Tomita-Takesaki modular automorphisms, Connes' Radon-Nikodym derivatives, Haagerup's standard form, or the Gelfand-Neumark-Segal construction, which allows us to define $\sigma^\mu$ and $(\mathrm{D}\mu : \mathrm{D}\nu)$ as explained above. Haagerup's standard form of $M$ can then be defined to be the space $\mathrm{L}_{1/2}(M)$, whereas the GNS construction of $M$ with respect to a weight $\mu$ is given by the map $x \mapsto x\mu^{1/2}$ defined whenever $x\mu^{1/2} \in \mathrm{L}_{1/2}(M)$, or, equivalently, $\mu(x^*x)$ is finite. To keep this paper self-contained we stick to the traditional approach.

6.9 The theory of representable functors assigns a precise meaning to the notion of an **I**-graded von Neumann algebra $\hat{M}$ generated by a family of generators and relations. The idea is to define a functor $F_M$ from **I**-graded von Neumann algebras to sets, prove its representability, and then use the Yoneda lemma to define $\hat{M}$ as the (unique) object that represents $F_M$. Thus $F_M$ will be isomorphic to the functor $\mathrm{Mor}(\hat{M}, -)$, i.e., for every **I**-graded von Neumann algebra $N$ the set $F(N)$ consists of all morphisms from $\hat{M}$ to $N$ and for every morphism $f \colon N \to O$ of **I**-graded von Neumann algebras the map of sets $F(f) \colon F(N) \to F(O)$ is given by the composition with $f$.

6.10 **Definition.** Given a von Neumann algebra $M$, the functor $F_M$ from the category of **I**-graded von Neumann algebras to the category of sets sends an **I**-graded von Neumann algebra $N$ to the set of all pairs $(f, g)$ such that $f$ is a morphism from $M$ to the 0-graded component of $N$, which is a von Neumann algebra, $g$ is a map that sends every $\mu \in \mathrm{W}_{\mathrm{sf}}(M)$ to a continuous morphism of groups $g_\mu \colon \mathbf{I} \to \mathbf{U}(f(p)Nf(p))$ ($p$ is the support of $\mu$) such that $g_\mu(a)$ has grading $a$ for all $a \in \mathbf{I}$, and finally the pair $(f, g)$ satisfies the two relations above concerning the modular automorphism group and the cocycle derivative. The functor $F_M$ sends a morphism $h \colon N \to O$ of **I**-graded von Neumann algebras to the map of sets $F_M(N) \to F_M(O)$ given by the composition of $f$ and $g$ with $h$.

6.11 **Theorem.** For any von Neumann algebra $M$ the functor $F_M$ defined above is representable.

6.12 *Proof.* The representable functor theorem (see Theorem 2.9.1 in Pareigis [11]) states that a functor $F$ from a complete category $C$ to the category of sets is representable if and only if $F$ is continuous (preserves small limits) and satisfies the solution set condition: there is a set $A$ of objects in $C$ such that for every object $X$ in $C$ and for every $x \in F(X)$ there are $W \in A$, $w \in F(W)$, and $h \colon W \to X$ such that $x = F(h)(w)$. We apply this theorem to the case when $C$ is the category of **I**-graded von Neumann algebras and $F$ is the functor $F_M$ constructed above.

6.13 The category of **I**-graded von Neumann algebras is complete (Guichardet [10] proves the result for ordinary von Neumann algebras, which immediately extends to the **I**-graded case). To prove that the functor $F_M$ preserves small limits it is sufficient to prove that $F_M$ preserves small products and equalizers. The functor $F$ preserves small products because a weight on a product decomposes into a family of weights on the individual factors, and the same is true for one-parameter groups of unitary elements. It preserves equalizers because the equalizer of two von Neumann algebras is their set-theoretical equalizer equipped with the restriction of the relevant structures. Finally, $F_M$ satisfies the solution set condition: for every **I**-graded von Neumann algebra $X$ and every element $x \in F_M(X)$ the **I**-graded von Neumann subalgebra $W$ of $X$ generated by the image of $M$ and all elements in the one-parameter families corresponding to all semifinite weights on $M$ is bounded in cardinality uniformly with respect to $X$, and all **I**-graded von Neumann algebras with cardinality at most some cardinal have a set of representatives of isomorphism classes. Thus the functor $F$ is representable. ∎

6.14 **Definition.** If $M$ is a von Neumann algebra, then the *core* of $M$ is the von Neumann algebra $\hat{M}$ that represents the functor $F_M$ defined above.

6.15 There are alternative ways to prove representability. For example, after the proof of continuity one can simply refer to the fact that the category of **I**-graded von Neumann algebras is locally presentable, as explained in a forthcoming paper by the author [47].

6.16 Alternatively, the universal property of crossed products (see Theorem 2 in Landstad [18]) allows us to prove that the crossed product of $M$ by the modular automorphism group of an arbitrary faithful semifinite weight on $M$ represents $F_M$. In particular, all of these crossed products for different weights are functorially isomorphic to each other and to the core, which can be seen directly by combining the universal property of crossed products with the Radon-Nikodym derivative of the corresponding weights. Hence the limit (or the colimit) of all such crossed products also represents the functor $F_M$ and gives us a functorial construction of the core that does not depend on a choice of a weight. This construction is similar to the first construction, except that here we get rid of the choice of a weight at a later stage.

6.17 Another way to prove representability is to construct an algebraically **I**-graded unital *-algebra $\check{M}$ generated by the above generators and relations and take its completion in the weakest topology that makes all of its representations continuous, as explained by Yamagami [34, the two paragraphs before Lemma 2.1] and Sherman [39, the paragraph after the formula (2.2)]. Here a representation of $\check{M}$ can be described as a morphism of algebraically **I**-graded unital *-algebras from $\check{M}$ to an arbitrary **I**-graded von Neumann algebra $N$ corresponding (in the obvious sense) to some element of $F(N)$. This construction is essentially an expansion of the usual proof of the representable functor theorem.

6.18 The algebra $\hat{M}$ has a canonical faithful semifinite operator-valued weight $\epsilon$ and a faithful semifinite trace $\tau$, which we describe briefly. Denote by $\epsilon$ the faithful semifinite operator-valued weight from $\hat{M}$ to $M$ corresponding to the embedding of $M$ into $\hat{M}$ and defined by the equality $\epsilon(x) = \int_\mathbf{R}(s \in \mathbf{R} \mapsto \theta_s(x) \in \hat{M}^+)$ for all $x \in \hat{M}^+$ and by $\tau$ the faithful semifinite trace on $\hat{M}$ defined by the equality $(\mathrm{D}(\mu \circ \epsilon) : \mathrm{D}\tau)_t = \mu^t$ for all $t \in \mathbf{I}$ and $\mu \in \mathrm{W}_{\mathrm{sf}}(M)$. We have $\tau \circ \theta_s = \exp(-s)\tau$. See Yamagami [34, §2], Falcone and Takesaki [38, Section 3], as well as Takesaki [41, §XII.6] for the relevant proofs.

6.19 We extend the construction of the core to a functor. First we need to define an appropriate notion of a morphism between von Neumann algebras. It turns out that in this situation a morphism from $M$ to $N$ is a pair $(f, T)$, where $f: M \to N$ is a usual morphism of von Neumann algebras and $T: \hat{\mathrm{L}}_0^+(N) \to \hat{\mathrm{L}}_0^+(M)$ is a faithful semifinite operator-valued weight associated to the morphism $f$, i.e., an ultraweakly lower semicontinuous (alias normal) additive homogeneous map such that $T(f(p)qf(p)^*) = pT(q)p^*$ for all $p \in M$ and $q \in \hat{\mathrm{L}}_0^+(N)$. Here $\hat{\mathrm{L}}_0^+(M)$ denotes the extended positive cone of $M$, which is a certain completion of $M^+$ that can be thought of as consisting of positive unbounded elements affiliated with $M$, as explained by Haagerup [16, Section 1] (he uses the notation $\hat{M}^+$ for $\hat{\mathrm{L}}_0^+(M)$, which conflicts with our notation for the positive part of the core of $M$). In the smooth case one can think of positive unbounded smooth functions on a smooth manifold. More generally, the positive cone $\mathrm{L}_a^+(M)$ of $\mathrm{L}_a(M)$ for $a \geq 0$ can be completed in a similar fashion to the extended positive cone $\hat{\mathrm{L}}_a^+(M)$ of $\mathrm{L}_a(M)$. For $a = 0$ we recover Haagerup's construction. For $a = 1$ the cone $\hat{\mathrm{L}}_1^+(M)$ consists of all weights on $M$ (not necessarily faithful or semifinite), in particular its semifinite part is precisely $\mathrm{W}_{\mathrm{sf}}(M)$. For other $a > 0$ the cone $\hat{\mathrm{L}}_a^+(M)$ is homeomorphic to $\hat{\mathrm{L}}_1^+(M)$, but the algebraic structures are different.

6.20 Semifiniteness and faithfulness for operator-valued weights are defined as for ordinary weights. Composition of such morphisms is componentwise and is well-defined because the composition of faithful respectively semifinite operator-valued weights is again faithful respectively semifinite. See Haagerup [16, 17] or Takesaki [41, §IX.4] for the relevant facts about operator-valued weights. Given a morphism $(f,T)\colon M \to N$ we define a morphism $\hat{f}\colon \hat{M} \to \hat{N}$ via the following formulas: $\hat{f}(x) = f(x)$ for all $x \in M$ and $\hat{f}(\mu^t) = (\mu \circ T)^t$ for all $\mu \in W_{\mathrm{sf}}(M)$ and $t \in \mathbf{I}$. Here $\mu \circ T \in W_{\mathrm{sf}}(N)$ because $\mu$ is a semifinite operator-valued weight from $M$ to $\mathbf{C}$ associated to the only morphism from $\mathbf{C}$ to $M$. The theory of operator-valued weights implies that such a mapping preserves all relations between generators (here we use the fact that $T$ is faithful) and hence defines a morphism $\hat{f}$ from $\hat{M}$ to $\hat{N}$ by the universal property of $\hat{M}$. Observe that $T$ has to be faithful because for a faithful semifinite weight $\mu$ on $M$ the one-parameter family $t \in \mathbf{I} \mapsto \mu^t \in \hat{M}$ takes value 1 at 0, hence its image under $(f,T)$, which is the one-parameter family $t \in \mathbf{I} \mapsto (\mu \circ T)^t \in \hat{N}$, must also take value 1 at 0, hence $\mu \circ T$ must be faithful, therefore $T$ must be faithful.

6.21 To obtain the analog of this construction in the smooth case we observe that morphisms of algebras correspond to smooth maps of manifolds going in the opposite direction. For simplicity we restrict to submersions, which have a good relative integration theory. An operator-valued weight then corresponds to a pushforward operation on (positive unbounded) smooth functions. Such a pushforward is given by integrating a given smooth function fiberwise with respect to a fixed (strictly) positive relative 1-density, i.e., a positive section of the line bundle of 1-densities associated to the relative tangent bundle. Composition can then be defined by observing that submersions are closed under composition and the relative tangent bundle of the composition of $f$ and $g$ can be computed as an extension of the pullback of the relative tangent bundle of $g$ along $f$ via the relative tangent bundle of $f$. The line bundle of relative 1-densities of $gf$ can then be computed as the tensor product of the line bundle of relative 1-densities of $f$ and the pullback of the line bundle of relative 1-densities of $g$ along $f$. Thus positive relative 1-densities $\mu$ and $\nu$ associated to $f$ and $g$ respectively can be composed by pulling back $\nu$ along $f$ and applying the above isomorphism. Observe that positive 1-densities on a smooth manifold $M$ are precisely the positive relative 1-densities associated to the submersion $M \to \mathrm{pt}$ and therefore can be composed with submersions $N \to M$ equipped with positive relative 1-densities yielding positive 1-densities on $N$, thus giving a morphism from the pullback along $N \to M$ of the line bundle of 1-densities over $M$ to the line bundle of 1-densities over $N$. This morphism gives a morphism between the line bundles of imaginary densities over $M \times \mathbf{I}$ and $N \times \mathbf{I}$, and the induced morphism between the corresponding algebras of sections $\hat{M}$ and $\hat{N}$ is the smooth analog of $\hat{f}$.

6.22 Morphisms from $M$ to $N$ in the above category can also be defined as morphisms $\hat{f}$ of $\mathbf{I}$-graded von Neumann algebras $\hat{f}\colon \hat{M} \to \hat{N}$. Such a morphism induces a morphism of its components in grading 0, i.e., a morphism $M \to N$ of von Neumann algebras. Furthermore, a semifinite weight $\mu$ on $M$ yields a one-parameter group of partial isometries $t \in \mathbf{I} \mapsto \mu^t \in \mathbf{I}$ whose value at $t \in \mathbf{I}$ has grading $t$. Under the map $\hat{f}$ such a group is mapped to a one-parameter group of partial isometries in $\hat{N}$ with the same property, which therefore necessarily comes from a semifinite weight $\nu$ on $N$. Given a element $y \in \hat{L}_0^+(N)$ we can define an element $z \in \hat{L}_0^+(M)$ by specifying the value of $\mu(z)$ on all semifinite weights $\mu$, which we define to be $\nu(y)$. The theory of operator-valued weights implies that the resulting map $T\colon \hat{L}_0^+(N) \to \hat{L}_0^+(M)$ is a faithful semifinite operator-valued weight associated to $f$ and the established correspondence from $\hat{f}$ to pairs $(f, T)$ is bijective. See Yamagami [34, the two paragraphs after Corollary 3.2] for the details. The smooth analog of this statement is easy to establish: once we have a map from 1-densities on $M$ to 1-densities on $N$, which can be thought of as a section of the tensor product of the pullback along $f\colon N \to M$ of the dual of the line bundle of 1-densities on $M$ and the line bundle of 1-densities on $N$, we can define a positive relative 1-density associated to $f$ by identifying the above tensor product with the line bundle of relative 1-densities of $f$.

6.23 Thus we obtain a functor from the category of von Neumann algebras and faithful semifinite operator-valued weights to the category of $\mathbf{I}$-graded von Neumann algebras, which turn out to be semifinite and are equipped with algebraic gadgets like $\tau$ and $\epsilon$ interacting in a certain way. The core functor is fully faithful, and hence is an equivalence of the domain category and its essential image. Thus the study of the category of arbitrary von Neumann algebras and faithful semifinite operator-valued weights reduces to the study of the category of certain $\mathbf{I}$-graded semifinite von Neumann algebras.

6.24 We define the space $\mathrm{L}_a(M)$ to be the $a$-graded component of $\hat{M}$ for all $a \in \mathbf{I}$. To define $\mathrm{L}_a(M)$ for $\Re a > 0$ we need to introduce unbounded elements. See Nelson [12, Section 2] or Terp [22, Chapter I] for the relevant definitions and proofs. Denote by $\bar{M}$ the completion of $\hat{M}$ in the $\tau$-measurable topology. All algebraic operations (including the grading) on $\hat{M}$ are continuous in this topology, hence $\bar{M}$ is a topological unital *-algebra. In fact, it is a **C**-graded *extended von Neumann algebra*, which means that its properties are very similar to those of von Neumann algebras (in particular, it admits a version of the Borel functional calculus), but it also has unbounded elements, which have infinite norm. We define positive elements of $\bar{M}$ in the same way as for von Neumann algebras. The set of all positive elements of $\bar{M}$ is the closure of $M^+$ in the $\tau$-measurable topology. The **I**-grading on $\bar{M}$ extends analytically to a **C**-grading. For every $a \in \mathbf{C}_{\Re \geq 0}$ we define $\mathrm{L}_a(M)$ as the $a$-graded component of $\bar{M}$. We require $a \in \mathbf{C}_{\Re \geq 0}$ because all other graded components are zero. See Terp [22, Chapter II] and Yamagami [34, §2] for details.

6.25 The construction $M \to \bar{M}$ is functorial if we restrict ourselves to the subcategory of von Neumann algebras whose morphisms are *bounded* faithful operator-valued weights. An operator-valued weight is *bounded* if it sends bounded elements to bounded elements. An element of $\hat{M}$ is called *bounded* if it belongs to the image of $M$ in $\hat{M}$. Alternatively, a bounded operator-valued weight $T$ associated to a morphism $f\colon M \to N$ of von Neumann algebras is simply a continuous positive morphism of $M$-$M$-bimodules $N \to M$, where the $M$-actions on $N$ come from $f$. Just as operator-valued weights generalize weights, bounded operator-valued weights generalize positive elements of the predual. The reason for this boundedness condition is that the induced map from $\mathrm{L}_1^+(M)$ to $\mathrm{L}_1^+(N)$ sends $\mu$ to $\mu \circ T$ for $\mu \in \mathrm{L}_1^+(M)$. This formula requires that $\mu \circ T \in \mathrm{L}_1^+(N)$, and if $\mu \circ T \in \mathrm{L}_1^+(N)$ for all $\mu \in \mathrm{L}_1^+(M)$, then $T$ is bounded. Morphisms in the restricted category from $M$ to $N$ can also be defined as morphisms of **C**-graded extended von Neumann algebras $\bar{M} \to \bar{N}$. See Yamagami [34, the second paragraph before Lemma 3.5] for the details. Thus we obtain a fully faithful functor from the category of von Neumann algebras and faithful bounded operator-valued weights to the category of **C**-graded extended von Neumann algebras.

6.26 The smooth analog of a bounded operator-valued weight is a *bounded* positive relative 1-density, which is defined as a positive relative 1-density whose fiberwise integral exists and is bounded as a smooth function on the base. We observe that such densities are closed under composition, hence they induce pullback maps between $a$-densities for $\Re a > 0$ in the same fashion as for $\Re a = 0$. Indeed, the smooth $\mathrm{L}_a$-space for $\Re a > 0$ is the space of integrable smooth $a$-densities, which for $a = 1$ consists precisely of bounded relative 1-densities associated to the map from the base to the point.

6.27 For all $a \in \mathbf{I}$ the space $\mathrm{L}_a(M)$ defined above is a functor from the category of von Neumann algebras and faithful semifinite operator-valued weights to the category of Banach spaces, which is the composition of the core functor and the functor that extracts the relevant graded component from an **I**-graded von Neumann algebra. Likewise, the space $\mathrm{L}_a(M)$ for $a \in \mathbf{C}_{\Re \geq 0}$ is a functor from the category of von Neumann algebras and faithful bounded operator-valued weights to the category of (quasi)Banach spaces.

6.28 As expected, $\mathrm{L}_0(M)$ and $\mathrm{L}_1(M)$ are naturally isomorphic to $M$ and $M_*$. The definition of $\mathrm{L}_a(M)$ implies that the multiplication on $\bar{M}$ induces a bilinear map $\mathrm{L}_a(M) \times \mathrm{L}_b(M) \to \mathrm{L}_{a+b}(M)$ for all $a$ and $b$ in $\mathbf{C}_{\Re \geq 0}$. The spaces $\mathrm{L}_a(M)$ together with these bilinear maps form an algebraically $\mathbf{C}_{\Re \geq 0}$-graded ring. In particular, every complex vector space $\mathrm{L}_a(M)$ is an $M$-$M$-bimodule. Moreover, the involution on $\bar{M}$ restricts to an antiisomorphism of $M$-$M$-bimodules $\mathrm{L}_a(M)$ and $\mathrm{L}_{\bar{a}}(M)$, hence the graded ring introduced above is a unital algebraically $\mathbf{C}_{\Re \geq 0}$-graded *-algebra. Since the multiplication on $\bar{M}$ is associative, we have a functorial map $m\colon \mathrm{L}_a(M) \otimes_M \mathrm{L}_b(M) \to \mathrm{L}_{a+b}(M)$. The first main theorem of this paper states that $m$ is an isomorphism of algebraic $M$-$M$-bimodules. Even though $\mathrm{L}_a(M)$ is not finitely generated as a right $M$-module, there is no need to complete the tensor product $\mathrm{L}_a(M) \otimes_M \mathrm{L}_b(M)$.

6.29 The second main theorem states that this isomorphism is an isometry. To make sense of this statement we introduce (quasi)norms on $\mathrm{L}_a(M)$ and on the tensor product. Suppose for a moment that $\Re a \leq 1$. Then there is a natural norm on $\mathrm{L}_a(M)$. If $\Re a = 0$, this norm is the restriction of the norm on $\hat{M}$. (In this case $\mathrm{L}_a(M)$ is a subset of $\hat{M}$.) If $\Re a > 0$, then the norm is given by the map $x \in \mathrm{L}_a(M) \mapsto \bigl((x^*x)^{1/2\Re a}(1)\bigr)^{\Re a}$. Here $x \in \mathrm{L}_a(M)$ and $x^* \in \mathrm{L}_{\bar{a}}(M)$, hence $x^*x \in \mathrm{L}_{2\Re a}^+(M)$. By $\mathrm{L}_d^+(M)$ we denote the intersection of $\bar{M}^+$ with $\mathrm{L}_d(M)$ for an arbitrary $d \in \mathbf{R}_{\geq 0}$. (If $d$ is not real, then the intersection is zero.) The Borel functional calculus extended to $\bar{M}$ and applied to the function $t \in \mathbf{R}_{\geq 0} \mapsto t^e \in \mathbf{R}_{\geq 0}$ yields a map $\mathrm{L}_d^+(M) \to \mathrm{L}_{de}^+(M)$ for all real $d \geq 0$ and $e \geq 0$, which is a bijection for $e > 0$. Thus $y = (x^*x)^{1/2\Re a} \in \mathrm{L}_1^+(M)$. Since $\mathrm{L}_1^+(M) = M_*^+$ (here $M_*^+$ is the set of all positive functionals in $M_*$), we have $y(1) \geq 0$ and $\|x\| = (y(1))^{\Re a} \geq 0$. This norm turns $\mathrm{L}_a(M)$ into a Banach space. If $\Re(a+b) \leq 1$, then we equip the tensor product $\mathrm{L}_a(M) \otimes_M \mathrm{L}_b(M)$ with the projective tensor norm and the second main theorem states that the map $m$ is an isometry of Banach spaces.

6.30 If $\Re a > 1$, then the same formula as above gives a quasinorm on $\mathrm{L}_a(M)$. Quasinorms satisfy a relaxed triangle inequality $\|x+y\| \le c(\|x\|+\|y\|)$ for some $c \ge 1$. In our case we can take $c = 2^{\Re a-1}$. The space $\mathrm{L}_a(M)$ is complete with respect to this quasinorm and therefore is a quasi-Banach space. If $\Re(a+b) > 1$, then we equip the tensor product $\mathrm{L}_a(M) \otimes_M \mathrm{L}_b(M)$ with the generalization of the projective tensor norm by Turpin [24] and the second main theorem states that the map $m$ is an isometry of quasi-Banach spaces.

6.31 The left multiplication map $\mathrm{L}_a(M) \to \mathrm{Hom}_M(\mathrm{L}_b(M), \mathrm{L}_{a+b}(M))$ that sends an element $x \in \mathrm{L}_a(M)$ to the map $y \in \mathrm{L}_b(M) \mapsto xy \in \mathrm{L}_{a+b}(M)$ is also an isomorphism of algebraic $M$-$M$-bimodules. Similarly to the tensor product, all elements in the algebraic internal hom above are automatically continuous in the norm topology, in particular they are bounded and can be equipped with the usual norm, which turns the above isomorphism into an isometry. This result was known before for the case of the continuous internal hom, only the automatic continuity part is new.

6.32 We recall some properties of the left and right supports and of the left and right polar decompositions of an element $x \in M$. The right support of $x$ is the unique projection that generates the left annihilator of the right annihilator of $x$. It can also be defined as the infimum of all projections $p$ such that $x = xp$. The left support of $x$ is the right support of $x$ in the opposite algebra of $M$. The involution exchanges the left and the right support. The left and right supports of $x$ are equal if $x$ is normal. In this case we refer to both of them as the support of $x$. The first definition of the right support of $x$ implies that $xy = 0$ if and only if $py = 0$, where $p$ is the right support of $x$. Likewise for the left support.

6.33 The right polar decomposition of $x$ is the unique pair $(y, z) \in M \times M$ such that $x = yz$, $y$ is a partial isometry, $z \ge 0$, and the right support of $y$ equals the right support of $x$. It follows that $y^*x = z = (x^*x)^{1/2}$, the left support of $y$ equals the left support of $x$, and the right support of $y$ equals the support of $z$. Recall that for any partial isometry $y$ its left support is equal to $yy^*$ and its right support is equal to $y^*y$. The left polar decomposition of $x$ is the right polar decomposition of $x$ in the opposite algebra of $M$. If $x = yz$ is the right polar decomposition of $x$, then $x = (yzy^*)y$ is the left polar decomposition of $x$. Thus the partial isometry parts of both polar decompositions coincide and we refer to both of them as the partial isometry part of $x$.

6.34 The notions of supports and polar decompositions extend to the elements of $\tilde{M}$, where $\tilde{M}$ is the completion of a von Neumann algebra $M$ with respect to the $\tau$-measurable topology for some faithful semifinite trace $\tau$. The definitions of all notions are the same. Recall that the bounded elements of $\tilde{M}$ (an element $z \in \tilde{M}$ is bounded if $z^*z \le d \cdot 1$ for some $d \in \mathbf{R}_{\ge 0}$) are precisely the elements in the image of $M$ in $\tilde{M}$. In particular, all projections and partial isometries in $\tilde{M}$ come from $M$. Hence the left and right supports of any element in $\tilde{M}$ belong to $M$. See the next section for the construction of polar decompositions.

6.35 For the case of the core of $M$ with its canonical trace we can say more about homogeneous elements, i.e., elements of $\mathrm{L}_a(M)$ for some $a \in \mathbf{C}_{\Re \ge 0}$. It turns out that if $w \in \mathrm{L}_a(M)$ for some $a \in \mathbf{C}_{\Re \ge 0}$, then both supports of $w$ are in $M$ (and not only in $\hat{M}$), the partial isometry part of $w$ is in $\mathrm{L}_{\Im a}(M)$ (recall that the imaginary part $\Im a$ of $a$ belongs to $\mathbf{I}$, the space of imaginary complex numbers), and the positive parts are in $\mathrm{L}^+_{\Re a}(M)$. See the next section for a proof of all these statements.

# 7 The algebraic tensor product and internal hom of $L_a$ and $L_b$.

7.1 First we need the following version of the Douglas lemma (see Theorems 1 and 2 in Douglas [9] for the original), which essentially explains how to perform division on elements of $L_a$-spaces. Various bits and pieces of this result are scattered throughout the literature, see the formula (2.1) and Lemma 4.2 in Junge and Sherman [43], Lemma 2.13 and Lemma 3.5 in Yamagami [34], Sections 1.2, 1.4, and 1.5 and Lemmas 2.2 and 2.4 in Schmitt [30], Lemma VII.1.6 (i) in Takesaki [41].

7.2 **Douglas lemma and polar decomposition for measurable operators.** For any von Neumann algebra $N$ with a faithful semifinite trace $\tau$ and any elements $x$ and $y$ in $\tilde{N}$, where $\tilde{N}$ is the completion of $N$ with respect to the $\tau$-measurable topology, there exists an element $p \in N$ such that $px = y$ if and only if there exists a $c \in \mathbf{R}_{\geq 0}$ such that $c^2 x^* x \geq y^* y$. Moreover, if the right support of $p$ is at most the left support of $x$, then such a $p$ is unique. In this case the norm of $p$ equals the smallest possible value of $c$, and the left support of $p$ equals the left support of $y$.

7.3 If $x$ and $y$ satisfy the stronger condition $x^* x = y^* y$, then $p$ is a partial isometry, the right support of $p$ equals the left support of $x$, and $p^* y = x$. In particular, all elements of $\tilde{N}$ have a unique right polar decomposition with the standard properties if we set $x = (y^* y)^{1/2}$.

7.4 *Proof.* To prove uniqueness, suppose that $px = qx$, where the right supports of $p \in N$ and $q \in N$ are at most the left support $z$ of $x$, i.e., $pz = p$ and $qz = q$. We have $(p - q)x = 0$, therefore $(p - q)z = 0$ and $p = pz = qz = q$. Uniqueness implies that the left support of $p$ is at most (hence equals) the left support of $y$. Indeed, if $z$ is the left support of $y$, then $px = y = zy = zpx$, thus $p = zp$.

7.5 To construct $p$, assume first that $x \geq 0$. Set $z_\epsilon = f_\epsilon(x)$, where $f_\epsilon(t)$ equals $t^{-1}$ for all $t \geq \epsilon$ and is 0 for all other $t$. Note that $z_\epsilon \in N$ for all $\epsilon > 0$ because $f_\epsilon$ is bounded. Furthermore, $q_\epsilon = xz_\epsilon = xf_\epsilon(x)$ is the spectral projection of $x$ to the set $[\epsilon, \infty)$. In particular, $q_\epsilon$ belongs to $N$ and has norm at most 1. Now $z_\epsilon^* y^* y z_\epsilon \leq z_\epsilon^* c^2 x^* x z_\epsilon = c^2 q_\epsilon^2 \leq c^2$, hence $yz_\epsilon \in N$ and its norm is at most $c$. Set $p = \lim_{\epsilon \to 0} yz_\epsilon$. Here the limit is taken over all $\epsilon > 0$ in the ultraweak topology. Recall that the unit ball of $N$ is compact in the ultraweak topology, therefore the limit exists, $px = y$, and $\|p\| \leq c$. By construction $pq = p$, where $q = \lim_{\epsilon \to 0}$ is the support of $x$, hence the condition on the right support of $p$ is satisfied.

7.6 For the case of general $x$ denote by $(r, u)$ the right polar decomposition of $x = ru$, where $r$ is a partial isometry in $N$ and $u \in \tilde{N}^+$. We use the above construction for $u \geq 0$ to find $q \in N$ such that $qu = y$, which is possible because $c^2 u^* u = c^2 x^* x \geq y^* y$. For $p = qr^*$ we have $p \in N$ and $px = qr^* x = qu = y$.

7.7 Fix the unique $p$ found above. We have $y^* y = x^* p^* px$. The inequality $c^2 x^* x \geq y^* y$ is true if and only if $x^*(c^2 - p^* p)x \geq 0$, equivalently $z^*(c^2 - p^* p)z \geq 0$, where $z$ is the left support of $x$. The latter inequality can be rewritten as $(cz)^2 \geq (pz)^*(pz) = p^* p$ and it is equivalent to $c^2 \geq p^* p$. Thus the minimum value of $c$ is precisely the norm of $p$.

7.8 Now suppose that $x^* x = y^* y$. We already proved that there is a unique $p \in N$ such that $px = y$ and the right support of $p$ is at most $q$, i.e., $pq = p$, where $q$ is the left support of $x$. Observe that $x^*(p^* p - 1)x = (px)^* px - x^* x = y^* y - y^* y = 0$, hence $q^*(p^* p - 1)q = 0$, which implies that $p^* p = (pq)^* pq = q^* q = q$, therefore $p^* p$ is a projection and $p$ is a partial isometry whose right support equals $q$ and $p^* y = p^* px = x$. ∎

7.9 **Corollary.** If $x \in L_a(M)$ and $y \in L_b(M)$ for some $a$ and $b$ in $\mathbf{C}_{\Re \geq 0}$ satisfy the conditions of the previous theorem for the core of $M$ with its canonical trace $\tau$, then $\Re a = \Re b$ (unless $y = 0$) and $p \in L_{b-a}(M)$. In particular, all elements of $L_a(M)$ have a unique right polar decomposition as the product of a partial isometry in $L_{\Im a}(M)$ and a positive element in $L_{\Re a}(M)$. Hence the left and right supports of an arbitrary element of $L_a(M)$ belong to $M$ because for an arbitrary partial isometry $u \in L_{\Im a}(M)$ the elements $u^* u$ and $uu^*$ belong to $M$.

7.10 *Proof.* We have $c^2 x^* x \geq y^* y$, where $c^2 x^* x \in L_{2\Re a}(M)$ and $y^* y \in L_{2\Re b}(M)$. The inequality is preserved under the action of $\theta_s$ for all $s \in \mathbf{R}$, hence $\Re a = \Re b$ unless $y = 0$. Suppose $px = y$ for some $p \in \hat{M}$. For all real $s$ we have $\theta_s(px) = \theta_s(p)\theta_s(x) = \theta_s(p)\exp(-sa)x$ and $\theta_s(px) = \theta_s(y) = \exp(-sb)y$, hence $\exp(sb - sa)\theta_s(p)x = y$. Combining this with the fact that $px = y$ we obtain $(p - \exp(sb - sa)\theta_s(p))x = 0$ and therefore $(p - \exp(sb - sa)\theta_s(p))e = 0$, where $e$ is the left support of $x$. Therefore $pe = \exp(sb - sa)\theta_s(p)e$. We have $pe = p$ and $\theta_s(p) = \theta_s(pe) = \theta_s(p)e$. Thus $p = \exp(sb - sa)\theta_s(p)$ and $\theta_s(p) = \exp(-s(b - a))p$, hence $p \in L_{b-a}(M)$. ∎

7.11 **Lemma.** If $M$ is a von Neumann algebra and $a \in \mathbf{C}_{\Re \geq 0}$, then any finitely generated left $M$-submodule $U$ of $\mathrm{L}_a(M)$ is generated by one element. The same is true if $U$ is countably generated and closed under countable sums. Likewise for right $M$-submodules.

7.12 *Proof.* Suppose that $U$ is generated by a family $u\colon I \to \mathrm{L}_a(M)$ for some finite set $I$. Set $x = (\sum_{i \in I} u_i^* u_i)^{1/2}$. Note that $u_i^* u_i \in \mathrm{L}_{2\Re a}^+(M)$ for all $i \in I$, therefore the sum is also in $\mathrm{L}_{2\Re a}^+(M)$, hence $x \in \mathrm{L}_{\Re a}^+(M)$. We have $x^* x = \sum_{i \in I} u_i^* u_i \geq u_i^* u_i$ for all $i \in I$. By the corollary above $u_i = q_i x$ for some $q_i \in \mathrm{L}_{\Im a}(M)$ for all $i \in I$. Choose an arbitrary faithful semifinite weight $\mu$ on $M$. Now $\mu^{-\Im a} \mu^{\Im a} = 1$ and therefore $u_i = q_i x = (q_i \mu^{-\Im a})(\mu^{\Im a} x)$. Thus all $u_i$ are left $M$-multiples of $y = \mu^{\Im a} x$. In the case of countable infinite families we must first rescale $u_i$ so that the sum in the definition of $x$ exists. The rest of the proof proceeds as above.

7.13 The proof will be complete when we show that $y \in U$. If $I = \emptyset$, then $y = 0 \in U$. Otherwise fix an element $k \in I$ and set $N = M \otimes \mathrm{End}(\mathbf{C}^I)$. Consider two elements $Y$ and $Z$ in $\mathrm{L}_a(N)$ such that $Y_{k,k} = y$, $Z_{i,k} = u_i$ for all $i \in I$ and all other entries of $Y$ and $Z$ are 0. Since $Y^* Y = Z^* Z$, there is a partial isometry $P \in N$ such that $Y = PZ$. In particular, $y = \sum_{i \in I} p_{k,i} u_i$, hence $y \in U$. In the countable infinite case we use the fact that $U$ is closed under countable sums. ∎

7.14 **Rank 1 theorem.** For any von Neumann algebra $M$, any right $M$-module $X$, any $a \in \mathbf{C}_{\Re \geq 0}$, and any element $z \in X \otimes_M \mathrm{L}_a(M)$ there exist $x \in X$ and $y \in \mathrm{L}_a(M)$ such that $z = x \otimes_M y$.

7.15 *Proof.* Represent $z$ as $\sum_{i \in I} u_i \otimes_M v_i$ for some finite set $I$ and some finite families $u\colon I \to X$ and $v\colon I \to \mathrm{L}_a(M)$. By the lemma above there exists an element $y \in \mathrm{L}_a(M)$ and a finite family $q\colon I \to M$ such that $v_i = q_i y$ for all $i \in I$. Now $z = \sum_{i \in I} u_i \otimes_M v_i = \sum_{i \in I} u_i \otimes_M q_i y = \sum_{i \in I} u_i q_i \otimes_M y = \left(\sum_{i \in I} u_i q_i\right) \otimes_M y$. Hence $x = \sum_{i \in I} u_i q_i$ and $y$ satisfy the requirements of the theorem. ∎

7.16 We are ready to prove the first main result of this paper.

7.17 **Algebraic tensor product isomorphism theorem.** For any von Neumann algebra $M$ and any $a$ and $b$ in $\mathbf{C}_{\Re \geq 0}$, the multiplication map $m\colon \mathrm{L}_a(M) \otimes_M \mathrm{L}_b(M) \to \mathrm{L}_{a+b}(M)$ is an isomorphism of algebraic $M$-$M$-bimodules. The inverse map is denoted by $n$ and is called the *comultiplication map*. Recall that $\otimes_M$ denotes the algebraic tensor product over $M$.

7.18 *Proof.* Since $m$ is a morphism of algebraic $M$-$M$-bimodules, it is enough to prove that $m$ is injective and surjective. To prove injectivity, suppose that we have an element $z \in \mathrm{L}_a(M) \otimes_M \mathrm{L}_b(M)$ such that $m(z) = 0$. By the rank 1 theorem there exist $x \in \mathrm{L}_a(M)$ and $y \in \mathrm{L}_b(M)$ such that $z = x \otimes_M y$. We want to prove that $x \otimes_M y = 0$. Note that $xy = m(x \otimes_M y) = m(z) = 0$. Since $xy = 0$, we have $py = 0$, where $p \in M$ is the right support of $x$. Thus $z = x \otimes_M y = xp \otimes_M y = x \otimes_M py = x \otimes_M 0 = 0$.

7.19 To prove surjectivity, suppose that we have an element $z \in \mathrm{L}_{a+b}(M)$. Choose $t \in \mathrm{L}_{\Im(a+b)}(M)$ and $h \in \mathrm{L}_1^+(M)$ such that $z = t h^{\Re(a+b)}$ is the right polar decomposition of $z$. Now $t h^{\Re a - \Im b} \in \mathrm{L}_a(M)$, $h^b \in \mathrm{L}_b(M)$, and $m(t h^{\Re a - \Im b} \otimes_M h^b) = z$. Thus the map $m$ is surjective and therefore bijective. The above formulas also give an explicit construction of the map $n$. ∎

7.20 To prove the second main result of this paper we recall some basic facts from the theory of tensor products of quasi-Banach spaces. If $X$ is a $p$-normed space for some $p \geq 1$, $Y$ is a $q$-normed space for some $q \geq 1$, and $r \geq 1$ is a real number, then we can introduce an $r$-seminorm on $X \otimes Y$: $u \in X \otimes Y \mapsto \|u\|_r := \sup_B \|B(u)\| \in \mathbf{R}_{\geq 0}$, where $B$ ranges over all linear maps from $X \otimes Y$ to some $r$-normed space $Z$ such that for all $x \in X$ and $y \in Y$ we have $\|B(x \otimes y)\| \leq \|x\| \cdot \|y\|$. A theorem by Turpin [24] states that for $r \geq p + q - 1$ this $r$-seminorm is an $r$-norm. In particular, for $p = q = r = 1$ we get the usual projective tensor norm.

7.21 **Algebraic tensor product isometry theorem.** In the notation of the previous theorem, equip the space $\mathrm{L}_a(M) \otimes_M \mathrm{L}_b(M)$ with the factor $r$-norm of the Turpin tensor $r$-norm on $\mathrm{L}_a(M) \otimes \mathrm{L}_b(M)$, where $r = \max(1, \Re(a+b))$. Then both $m$ and $n$ are isometries.

7.22 *Proof.* It is enough to prove that both $m$ and $n$ are contractive. To prove that $n$ is contractive, suppose that $z \in \mathrm{L}_{a+b}(M)$. In the notation of the previous proof we have $\|n(z)\| = \|t h^{\Re a - \Im b} \otimes_M h^b\| \leq \|t h^{\Re a - \Im b}\| \cdot \|h^b\| =$ $\|h^{\Re a}\| \cdot \|h^{\Re b}\| = h(1)^{\Re a} \cdot h(1)^{\Re b} = h(1)^{\Re(a+b)} = \|h^{\Re(a+b)}\| = \|th^{a+b}\| = \|z\|$. Note that $\|th^{\Re a - \Im b}\| = \|th^{\Re a}\| = \|h^{\Re a}\|$ because $t$ is a partial isometry in $M$ such that $t^*t$ equals the support of $h$.

7.23 The multiplication map $\mathrm{L}_a(M) \otimes \mathrm{L}_b(M) \to \mathrm{L}_{a+b}(M)$ satisfies the condition on $B$ in the definition of the Turpin $r$-norm because of the Hölder-Kosaki inequality (Kosaki [27, Lemma 1]): $\|xy\| \leq \|x\| \cdot \|y\|$ for all $x \in \mathrm{L}_a(M)$ and $y \in \mathrm{L}_b(M)$. Hence the above map is contractive and therefore $m$ is also contractive. ∎

7.24 **Remark.** The (quasi)norm on $\mathrm{L}_a(M)$ determines $\Re a$ unless $M = \mathbf{C}$ or $M = 0$. For $\Re a \leq 1$, this follows from Corollary 6.5 in Kosaki [29], for $\Re a > 1$ this follows from the Aoki-Rolewicz theorem (see Section 2 in Kalton [42]). Thus in the generic case the number $r$ in the statement of the above theorem is determined by the Banach structures on $\mathrm{L}_a(M)$ and $\mathrm{L}_b(M)$. In the two special cases the norm on $\mathrm{L}_{a+b}(M)$ coincides with the norms on $\mathrm{L}_a(M)$ and $\mathrm{L}_b(M)$. Thus it is always possible to define the quasinorm in the above theorem using just the (quasi)norms without any reference to $a$ and $b$.

7.25 We now prove the corresponding result for algebraic homomorphisms. The only novelty of the theorem below is that we remove the boundedness condition on homomorphisms. Everything else has been proved before, see Theorem 2.5 in Junge and Sherman [43], Proposition 2.10 in Yamagami [34], and Proposition II.35 in Terp [22]. We start by proving an automatic continuity lemma, which then implies the desired result.

7.26 **Automatic continuity lemma.** Suppose $M$ is a von Neumann algebra, $a \in \mathbf{C}_{\Re \geq 0}$, and $X$ is a topological vector space equipped with a structure of a right $M$-module such that for any $x \in X$ the map $p \in M \mapsto xp \in X$ is continuous in the norm topology on $M$. Then every morphism $T$ of algebraic right $M$-modules from $\mathrm{L}_a(M)$ to $X$ is continuous if $\mathrm{L}_a(M)$ is equipped with the (quasi)norm topology.

7.27 *Proof.* Consider any family $u: I \to \mathrm{L}_a(M)$ such that the sum $\sum_{i \in I} \|u_i\|^2$ exists. We will now prove that the map $T$ sends $u$ to an arbitrarily small neighborhood $V$ of 0 in $X$ if we rescale $u$ by any $a \in [0, \epsilon]$ for some fixed $\epsilon > 0$. Set $v = \left( \sum_{i \in I} u_i u_i^* \right)^{1/2} r$, where $r$ is a fixed arbitrary unitary element of $\mathrm{L}_{\Im a}(M)$. The sum converges in the norm topology by the definition of $u$. Since $vv^* \geq u_i u_i^*$ for all $i \in I$, the Douglas lemma gives us a family $p: I \to M$ such that for all $i \in I$ we have $u_i = vp_i$ and $\|p_i\| \leq 1$. Now $T(u_i) = T(vp_i) = T(v)p_i$. The map $q \in M \mapsto T(v)q \in X$ is continuous by the definition of $T$, therefore the preimage of any neighborhood $V$ of 0 in $X$ must contain a ball of some radius $\epsilon > 0$ centered at the origin, in particular $T(au_i) = T(avp_i) = T(v)ap_i \in V$ for all $a \in [0, \epsilon]$ and all $i \in I$ because $\|ap_i\| \leq \epsilon$. This proves the above claim.

7.28 If $T$ is not continuous, then there is at least one neighborhood $V$ of $0 \in X$ such that its preimage under $T$ does not contain a ball of a positive radius. In particular, for any $\epsilon > 0$ we can find $w_\epsilon \in \mathrm{L}_a(M)$ such that $\|w_\epsilon\| \leq \epsilon^2$ and $T(w_\epsilon) \notin V$. Fix such a neighborhood $V$ and a sequence $\epsilon: I \to (0, \infty)$ such that $\sum_{i \in I} \epsilon_i^2$ exists. We construct a sequence $u: I \to \mathrm{L}_a(M)$ such that $\sum_{i \in I} \|u_i\|^2$ exists and $T(\|u_i\| \cdot u_i) \notin V$ by setting $u_i = \|w_{\epsilon_i}\|^{-1/2} w_{\epsilon_i}$. Indeed, $\|u_i\| = \|w_{\epsilon_i}\|^{-1/2} \|w_{\epsilon_i}\| = \|w_{\epsilon_i}\|^{1/2} \leq \epsilon_i$, hence $\sum_{i \in I} \|u_i\|^2 \leq \sum_{i \in I} \epsilon_i^2$ exists and $T(\|u_i\| \cdot u_i) = T(\|w_{\epsilon_i}\|^{1/2} \|w_{\epsilon_i}\|^{-1/2} w_{\epsilon_i}) = T(w_{\epsilon_i}) \notin V$. This contradicts the statement proved in the previous paragraph because $a = \|u_i\|$ can become smaller than any fixed $\epsilon > 0$ and $T(au_i) = T(w_{\epsilon_i})$ must be in $V$ for such $i$. ∎

7.29 **Algebraic internal hom isomorphism theorem.** For any von Neumann algebra $M$ and any $a$ and $b$ in $\mathbf{C}_{\Re \geq 0}$ the left multiplication map $\mathrm{L}_a(M) \to \mathrm{Hom}_M(\mathrm{L}_b(M), \mathrm{L}_{a+b}(M))$ $(x \mapsto (y \mapsto xy))$ is an isomorphism of algebraic $M$-$M$-bimodules. The inverse map is called the *left comultiplication map*. Recall that $\mathrm{Hom}_M$ denotes the algebraic internal hom of right $M$-modules, the left $M$-action on the internal hom is induced by the left $M$-action on $\mathrm{L}_{a+b}(M)$, and the right $M$-action is induced by the left $M$-action on $\mathrm{L}_b(M)$.

7.30 *Proof.* If for some $x \in \mathrm{L}_a(M)$ we have $xy = 0$ for all $y \in \mathrm{L}_b(M)$, then $xp = 0$ for all $\sigma$-finite projections $p \in M$. Thus $x = 0$ and the map is injective.

7.31 To establish surjectivity, observe that the space $\mathrm{L}_{a+b}(M)$ satisfies the conditions of the automatic continuity lemma because $\|xm\|^2 = \|xmm^*x^*\| \leq \|m\|^2 \|x^*x\| = \|m\|^2 \|x\|^2$ for all $m \in M$ and $x \in \mathrm{L}_{a+b}(M)$. Hence all elements of $\mathrm{Hom}_M(\mathrm{L}_b(M), \mathrm{L}_{a+b}(M))$ are continuous in the (quasi)norm topology, i.e., bounded. We now invoke Theorem 2.5 in Junge and Sherman [43], which states that every bounded homomorphism of right $M$-modules from $\mathrm{L}_b(M)$ to $\mathrm{L}_{a+b}(M)$ is the left multiplication by an element of $\mathrm{L}_a(M)$. ∎

7.32 **Algebraic internal hom isometry theorem.** In the notation of the previous theorem, equip the space $\mathrm{Hom}_M(\mathrm{L}_b(M), \mathrm{L}_{a+b}(M))$ with the quasinorm $\|f\| = \sup_{\|y\| \leq 1} \|f(y)\|$, where $y \in \mathrm{L}_b(M)$. Then both $m$ and $n$ are isometries.

7.33 *Proof.* For $x \in \mathrm{L}_a(M)$ we have $\|m(x)\| = \sup_{\|y\| \leq 1} \|xy\| \leq \sup_{\|y\| \leq 1} \|x\| \cdot \|y\| \leq \|x\|$, thus $\|m(x)\| \leq \|x\|$. It remains to prove that $\|x\| \leq \|m(x)\|$. We can assume that $x \neq 0$.

7.34 If $\Re a \neq 0$, then we construct $y \in \mathrm{L}_b(M)$ such that $\|y\| \neq 0$ and $\|xy\| = \|x\| \cdot \|y\|$. Set $z = (x^*x)^{1/2\Re a}$ and $y = z^b$. Note that $x = uz^a$ for some partial isometry $u \in M$. Now $\|xy\| = \|uz^ay\| = \|z^ay\| = \|z^az^b\| = \|z^{a+b}\| = \|z\|^{\Re a + \Re b} = \|z\|^{\Re a}\|z\|^{\Re b} = \|z^a\| \cdot \|z^b\| = \|uz^a\| \cdot \|y\| = \|x\| \cdot \|y\|$.

7.35 In the case $\Re a = 0$ the above strategy does not work because such an element $y$ might not exist. Instead for every $c \in [0, \|x\|)$ we construct a nonzero $y \in \mathrm{L}_b(M)$ such that $\|xy\| \geq c\|y\|$. Set $z = (xx^*)^{1/2}$. Note that $x = zu$ for some partial isometry $u \in \mathrm{L}_a(M)$. Denote by $p$ the spectral projection of $z$ for the set $[c, \infty)$ so that $zp \geq c$. Choose a nonzero $w \in \mathrm{L}_1^+(M)$ with the support at most $p$ and set $y = u^*w^{a+b}$. Now $\|xy\| = \|zuu^*w^{a+b}\| = \|zw^{a+b}\| = \|zpw^{a+b}\| = \|w^{a+b}pzzpw^{a+b}\|^{1/2} \geq \|w^{a+b}c^2w^{a+b}\|^{1/2} = c\|w^{a+b}\| = c\|u^*w^{a+b}\| = c\|y\|$. ∎

7.36 **Remark.** For any von Neumann algebra $M$ the spaces $\mathrm{L}_a(M)$ for all $a \in \mathbf{C}_{\Re \geq 0}$ can be organized into a smooth bundle L of quasi-Banach $M$-$M$-bimodules over $\mathbf{C}_{\Re \geq 0}$ by postulating that a section of this bundle is smooth if it is locally of the form $a \in \mathbf{C}_{\Re \geq 0} \mapsto f(a)\mu^a$ for some smooth function $f \colon \mathbf{C}_{\Re \geq 0} \to M$ and some $\mu \in \mathrm{L}_1^+(M)$. The tensor product isomorphisms can be combined into a smooth bundle isomorphism $i^*\mathrm{L} \otimes_M j^*\mathrm{L} \to k^*\mathrm{L}$, where $i, j, k \colon \mathbf{C}_{\Re \geq 0} \times \mathbf{C}_{\Re \geq 0} \to \mathbf{C}_{\Re \geq 0}$ are respectively projections on the first and the second component and the addition map. Likewise, the internal hom isomorphisms can be combined into a smooth isomorphism $i^*\mathrm{L} \to \mathrm{Hom}_M(j^*\mathrm{L}, k^*\mathrm{L})$. In particular, if we restrict the tensor product isomorphism from $\mathbf{C}_{\Re \geq 0}$ to $\mathbf{I}$ and allow $\mu \in \mathrm{W}_{\mathrm{sf}}(M)$ in the definition of smooth sections, we obtain a convolution product on the space of distributional sections of L restricted to $\mathbf{I}$ with bounded Fourier transform, which turns it into an $\mathbf{I}$-graded von Neumann algebra (the grading is the composition of the isomorphism $\mathbf{R} \to \mathrm{Hom}(\mathbf{I}, \mathbf{U})$ and the multiplication action), which is canonically isomorphic to $\hat{M}$. This resembles the approach used by Falcone and Takesaki [38, Theorem 2.4] to construct $\hat{M}$. The case of $\bar{M}$ is similar if we do not restrict to $\mathbf{I}$.

## 8 Future directions and applications.

8.1 We briefly indicate how results in this paper can be extended to the setting of $\mathrm{L}_d(M)$-modules introduced by Junge and Sherman [43]. Details are to be found in a forthcoming paper by the author [46]. $\mathrm{L}_d(M)$-modules for a von Neumann algebra $M$ and $d \in \mathbf{R}_{\geq 0}$ generalize Hilbert W*-modules for $d = 0$ and representations of $M$ on Hilbert spaces for $d = 1/2$. They can be defined as algebraic $M$-modules equipped with an $\mathrm{L}_{2d}(M)$-valued inner product satisfying a natural set of properties. For any fixed $d \in \mathbf{R}_{\geq 0}$ the category of $\mathrm{L}_d(M)$-modules is a W*-category and the categories for different values of $d$ are all equivalent to each other. This is implicitly stated by Junge and Sherman in [43], see Theorem 3.6 and Proposition 6.1. Their proof is analytical in nature and occupies the entire Section 4.

8.2 In [46] we give an alternative algebraic proof of this result, which is very similar in spirit to this paper. The equivalence from $\mathrm{L}_d(M)$-modules to $\mathrm{L}_{d+e}(M)$-modules is given by the algebraic tensor product with $\mathrm{L}_e(M)$ as an $M$-$M$-bimodule. The inverse functor is given by the algebraic internal hom from $\mathrm{L}_e(M)$. The space $\mathrm{L}_d(M)$ has a canonical structure of an $\mathrm{L}_d(M)$-module, which allows us to recover the main results of this paper as a particular case, at least if $a$ and $b$ are real. The proof for $\mathrm{L}_d(M)$-modules, however, needs this particular case as a prerequisite.

## 9 Acknowledgments.

I thank Boris Ettinger for several suggestions about noncommutative $L_a$-spaces, David Sherman for sending me a copy of Terp's unpublished preprints [22] as well as many other papers, Joachim Cuntz for pointing out a missing part in the proof of the lemma about finitely generated left submodules of $L_a(M)$, André Henriques for numerous fruitful discussions about von Neumann algebras, and last and most importantly, my doctoral advisor Peter Teichner for introducing me to this area, suggesting to explore some conjectures proved in this paper, conducting numerous invaluable discussions, and supporting me through the course of my research at the University of California, Berkeley. This work was partially supported by the SFB 878 grant.

## 10 References.